\documentclass{amsart}
\usepackage{bbm}
\usepackage{dsfont}
\usepackage{relsize}
\usepackage{mathtools}
\usepackage{amsthm}
\usepackage[breaklinks]{hyperref}
\usepackage{blindtext}
\usepackage{amsmath,amssymb}
\usepackage{autonum}
\usepackage{enumitem}
\usepackage[usenames,dvipsnames,svgnames]{xcolor}
\usepackage{mathrsfs}
\usepackage{csquotes}
\usepackage[]{hyperref}
\hypersetup{colorlinks=true,linkcolor=WildStrawberry,citecolor=cyan}
\newtheorem{thm}{Theorem}[section]
\newtheorem{cor}[thm]{Corollary}
\newtheorem{lem}[thm]{Lemma}

\theoremstyle{definition}
\newtheorem{defn}[thm]{Definition}
\theoremstyle{remark}
\newtheorem{rem}[thm]{Remark}
\theoremstyle{example}

\numberwithin{equation}{section}

\newcommand{\To}{\longrightarrow}

\renewcommand{\a}{\alpha}
\renewcommand{\b}{\beta}

\makeatletter
\renewcommand\paragraph{\@startsection{paragraph}{4}{\z@}%
	{-2.5ex\@plus -1ex \@minus -.25ex}%
	{1.25ex \@plus .25ex}%
	{\normalfont\normalsize\centering\bfseries}}
\makeatother
\begin{document}
	\title[Mixed local and nonlocal elliptic equation with gradient degeneracy]
	{REGULARITY OF SOLUTIONS TO VARIABLE-EXPONENT DEGENERATE MIXED FULLY NONLINEAR LOCAL AND NONLOCAL EQUATIONS}
	\author[P.\,Oza,\,\, J.\,Tyagi]
	{ Priyank Oza,\,\,Jagmohan Tyagi }
	\address{Priyank\,Oza \hfill\break
		Indian Institute of Technology Gandhinagar \newline
		Palaj, Gandhinagar Gujarat, India-382055.}
	\email{priyank.k@iitgn.ac.in, priyank.oza3@gmail.com}
	\address{Jagmohan\,Tyagi \hfill\break
		Indian Institute of Technology Gandhinagar \newline
		Palaj, Gandhinagar Gujarat, India-382055.}
	\email{jtyagi@iitgn.ac.in, jtyagi1@gmail.com}
	\thanks{Submitted \today.  Published-----.}
\begin{abstract}
		We consider a class of variable-exponent mixed fully nonlinear local and nonlocal degenerate elliptic equations, which degenerate along the set of critical points, $C:=\big\{x:\,Du(x)=0\big\}.$ Under general conditions, first, we establish the Lipschitz regularity of solutions using the Ishii-Lions viscosity method when the order of the fractional Laplacian, $s\in\big({1}/{2},1\big)$ (Theorem \ref{Main 1}). Due to the inapplicability of the comparison principle for the equations under consideration, one can not employ the classical Perron's method for the existence of a solution. Utilizing the Lipschitz estimates established in Theorem \ref{Main 1} and \enquote{vanishing viscosity} method, we prove the existence of a solution. We further prove interior $C^{1,\delta}$ regularity of viscosity solutions using an improvement of the flatness technique when $s$ is close enough to $1$ (Theorem \ref{main 2}). 
	\end{abstract}
	\maketitle
	
	{\hypersetup{linkcolor=black}
	\tableofcontents}
	\section{Introduction} We establish interior $C^{1,\delta}$ regularity of viscosity solutions to the following mixed local and nonlocal degenerate equations:
	\begin{align}\label{eq 0.1}
		|Du|^{\gamma(x)}\big(\a F\big(D^2u\big)-\b(-\Delta)^su\big)+b\cdot Du|Du|^{\gamma(x)}=f &\text{ in } \,{B_1},
	\end{align}
	where $\a$ and $\b$ are non-negative constants, $B_1$ is the unit ball in $\mathbb{R}^N$ centred at the origin, $\gamma\in C(\overline{B}_1,\mathbb{R})$ with ${\inf_{x\in B_1}}\gamma(x)\geq 0,$ $f\in L^\infty({B}_1,\mathbb{R})$,\,$Du$ denotes the gradient of $u,$ $F$ is a uniformly elliptic operator (fully nonlinear), $(-\Delta)^s$ denotes the fractional Laplacian with order $s\in (0,1)$ and $b:B_1\To\mathbb{R}^N$ is a bounded Lipschitz continuous function which is specified later. 
	
Results concerning the regularity of solutions to PDEs have gained much attention in the last several decades. Equations of the form \eqref{eq 0.1} have been investigated widely in the last decade. We begin by mentioning the pioneering works of I. Birindelli and F. Demengel in a series of papers for degenerate equations pertaining to fully nonlinear uniformly elliptic local operators: comparison principle and Liouville-type results \cite{Birin 1}, regularity and uniqueness of eigenvalues and eigenfunctions \cite{Birin 3, Birin 2}, $C^{1,\alpha}$ regularity in the radial case \cite{Birin 4}. One may also see \cite{Davila, Davila 2, Imbert 1} for ABP estimates and Harnack's inequalities for degenerate, singular fully nonlinear equations. It is worth mentioning that when $\gamma(x)\equiv 0,$ $\beta=0$ and $b\equiv 0$ in \eqref{eq 0.1}, the regularity results are well known, for instance, see \cite{Cabre}.

	%
	
When $\beta=0,\,b\equiv 0$ in $B_1,$ $\gamma(x)\equiv\gamma\geq0$ is a constant, \eqref{eq 0.1} reads as follows:
	\begin{align}\label{a=1}
		|Du|^\gamma F\big(D^2u\big)=f &\text{ in } \,{B_1}.
	\end{align}
	This equation was investigated by Imbert and Silvestre \cite{Imbert}. In this paper, they proved the interior regularity of solutions for a class of equations given by \eqref{a=1}. One may also see \cite{Arau, Birin Leo}, which are closely related to \cite{Imbert}. Ara\'ujo et al. \cite{Arau} established the sharp regularity of solutions to \eqref{a=1} when $\gamma(x)\equiv\gamma>0$ is a constant. We refer to the works of Teixeira \cite{Teixeira 1, Teixeira 2, Teixeira 3} for related regularity results in a different context. Birindelli et al. \cite{Birin Leo} established $C^{1,\delta}$ regularity of fully nonlinear degenerate equations with superlinear and sub-quadratic Hamiltonian terms. More precisely, authors considered equations of the form
	\begin{align}
		-|Du|^\gamma F(D^2u)+b|Du|^\beta=f \text{ in }\Omega,
	\end{align} 
	where $b, f$ are continuous functions in $\Omega\subset \mathbb{R}^N$ with $C^2$-boundary, $\gamma>-1$ and $\gamma+1<\beta\leq \gamma+2.$

	There are several other works, which are devoted to the degenerate elliptic equations with gradient degeneracy. We first mention a few of them, which are closely related to this paper.

	\begin{enumerate}
	\item[(i)] When $\b=0,\,b\equiv0,\,\gamma(x)\equiv\gamma\geq0$ is a constant, Imbert and Silvestre \cite{Imbert} established interior $C^{1,\delta}$ regularity of solutions to \eqref{eq 0.1}.
	\item[(ii)] For $\a=0,$ $b\equiv 0, \gamma(x)\equiv\gamma\geq0$ is a constant, dos Prazeres and Topp \cite{Topp} investigated interior $C^{1,\delta}$ regularity of solutions to \eqref{eq 0.1}, where the nonlocal term is nonlinear.
	\item[(iii)] Bronzi et al. \cite{Bronzi} studied \eqref{eq 0.1}, where $\gamma(x)$ is a continuous function, $\b=0$ and $b\equiv 0$. They established interior $C^{1,\delta}$ regularity. Inspired by the results of Teixeira \cite{Teixeira H, Teixeira 1}, they also established sharp geometric estimates.
	\item[(iv)] Birindelli and Demengel \cite{C1beta} discussed \eqref{eq 0.1} in the case when $\b=0,\gamma(x)\equiv\gamma\geq0$ is a constant, i.e., fully nonlinear local equation with a first-order term.
	\end{enumerate} 
	
	Very recently, De Filippis \cite{Filippis} introduced degeneracies of double phase type to fully nonlinear equation
	\begin{align}\label{Filippis}
		\left(|Du|^{\gamma_1}+a(x)|Du|^{\gamma_2}\right)F\big(D^2u\big)=f(x) \text{ in } \Omega,	
	\end{align}
	for some non-negative continuous function $a,$ $\gamma_2\geq \gamma_1>0$ and bounded domains $\Omega$ in $\mathbb{R}^N,$ $N\geq 2.$ De Filippis proved $C^{1,\delta}$ local regularity of viscosity solutions of \eqref{Filippis}. In another paper, De Filippis \cite{Filippis 1} studied the existence and regularity of solutions to free transmission problems, see \cite{Filippis 1}. One may also see \cite{Huaroto, Jesus} for the related works. In particular, Huaroto et al. \cite{Huaroto} considered the degenerate free transmission problem:
	\begin{align}
		|Du|^{\theta_1\mathbbm{1}_{\{u>0\}}+\theta_2\mathbbm{1}_{\{u<0\}}}F(D^2u)=f(x) \text{ in }\Omega \text{ with }\theta_1,\,\theta_2> 0.	
	\end{align} 
	Jesus \cite{Jesus} studied the following free transmission problem:
	\begin{align}
		|Du|^{\beta(x,u,Du)}F(D^2u)=f(x) \text{ in }B_1 \text{ with }\beta\geq 0.
	\end{align}
We also mention the works of da Silva et al. \cite{da Silva} and Fang et al. \cite{Fang}, who established the local $C^{1,\delta}$ regularity of solutions to fully nonlinear local equations with gradient degeneracy of the form 
	\begin{align}
		|Du|^{\gamma_1(x)}+a(x)|Du|^{\gamma_2(x)} \text{ with } \gamma_2(x)\geq \gamma_1(x)\geq 0 \text{ in }\Omega\subset\mathbb{R}^N,\,N\geq 2.
	\end{align}

The variational counterpart of equations of the form \eqref{Filippis} has been widely studied by several authors. The pioneering work is due to Marcellini \cite{var1, var2}. We refer to \cite{Colombo 1, Colombo 2, Colombo 3} for the works on double-phase energy models. More precisely, Colombo et al. \cite{Colombo 2} established the $C^{1,\delta}$ regularity of the bounded minimisers of the functional
\begin{align}\label{model var}
	u\mapsto\int(|Dw|^p+a(x)|Dw|^q)dx,
\end{align} 
where $a(x)$ is a H\"older continuous function and $1<p<q.$ Further, the same authors \cite{Colombo 3} studied the double phase variational problems having the model form \eqref{model var}.
Recently, Baroni et al. \cite{Colombo 1} proved the sharp regularity for a class of functionals with double phases. We recall that the operator $$G(Du, D^2u)\coloneqq|Du|^\gamma F(D^2u)$$ can be viewed as a non-variational extension of the $p$-Laplacian, see \cite{Attouchi, Banerjee}, which occurs in the theory of stochastic games and has applications to image processing. From an application point of view, degenerate equations with variable-exponent are quite important. These equations appear in the study of image enhancement. One may see, for instance, Chen et al. \cite{image en} and the references therein for the details. They proposed an alternative model for image restoration allowing the exponent to be a variable instead of a fixed constant. The variational approach to the image restoration problem consists of minimizing the functional
	\begin{align}
		\int_{B_1}\bigg(|Du|^{\gamma(x)}+\frac{\sigma}{2}|u-I|^2\bigg)dx,
	\end{align}
	where $u\colon B_1\subset\mathbb{R}^2\To\mathbb{R}$ is the true image and $I\colon B_1\subset\mathbb{R}^2\To\mathbb{R}$ is the observed image. Here, $I\coloneqq u(x)+\xi,$ with $\xi$ being a noise term. For the different choices of $\gamma,$ distinct algorithms need to be followed, see \cite{image en} for the details. For instance, the method for the case when $\gamma\equiv 1,$ is known as total variation (TV) based regularization, which was proposed by Rudin et al. \cite{ROF}. Each of these algorithms has its drawbacks. More precisely, when $\gamma(x)\equiv 1,$ TV-based solutions are piecewise constant. Even though $\gamma(x)\equiv 2,$ resolves this issue, the edges are not preserved. On the other hand, taking fixed $1< \gamma < 2$ may destroy the edges, see \cite{Obliterate} for the details. This elaborates on the importance of the study of variable-exponent degenerate equations. These equations also have applications in the study of electrorheological fluids \cite{electro 1, electro 2, electro 3} as well as in the thermistor problem \cite{thermo}. 
	
	It is important to note that the presence of fully nonlinear operator in the equations under consideration restricts the use of \enquote{linearization} techniques and bootstrap arguments to get regularity results. This has been exploited in several papers. For a brief overview of the bootstrap technique, one may see the book by Ambrosetti and Malchiodi \cite{Boot}, in particular, Theorem 1.16 therein. We refer to \cite{Felmer 18, Felmer 19, Oza Hopf, Pozza, Quaas, Quaas 2, Tyagi 3, Tyagi 1, Tyagi 2} for the existence and qualitative questions pertaining to extremal Pucci's equations. One may also see \cite{Tyagi degenerate} for fully nonlinear degenerate equations pertaining to partial trace operators, $\mathcal{P}_k^{\pm}$ with sublinear gradient term.

	In a pioneering work in the direction of establishing $C^{1,\delta}$ interior regularity of degenerate fully nonlinear equations, i.e., when $\beta=0,b\equiv 0$ and $\gamma(x)\equiv\gamma\geq0$ is a constant in \eqref{eq 0.1}, we refer to the work of Imbert and Silvestre \cite{Imbert}. In their proof, Lemma \ref{Lem Imbert} (see next) plays a major role. 
	
We mention that even in a very particular case, i.e., $\gamma(x)\equiv\gamma$ a positive constant, $\a=0$ and $b\equiv 0,$ the comparison principle does not hold for \eqref{eq 0.1}. One may see \cite{Topp} for a counterexample. Therefore one can not use the classical Perron's method to establish the existence of a viscosity solution. 
	
Recently, the elliptic equations involving mixed local (linear) and non-local operators without gradient degeneracy have gained much attention, see \cite{Biagi, Biagi 1, Biagi 2, Chen, Chen 2} and \cite{Anup mixed} for mixed local (fully nonlinear) and nonlocal equations. We mention the work \cite{Oza Heisenberg}, where we established the existence of viscosity solutions to mixed fully nonlinear local and nonlocal equations in the Heisenberg group using the Perron's method. We further proved the H\"older regularity of the viscosity solutions. These operators occur naturally in the study of plasma physics \cite{Plasma}, population dynamics \cite{Dipierro} and many more. Mou \cite{Mou} established the H\"older regularity of solutions to second order non-translation invariant integro-PDE. For the regularity results for mixed local and nonlocal $p$-Laplace equations, we refer to \cite{Biagi p, Prashanta 2, Prashanta 1, Mingione} and the references therein. In particular, Mingione \cite{Mingione} established the maximal regularity of solutions to mixed local and nonlocal operators of $p$-Laplacian type. In addition to the above works, we refer to a survey paper by Ros-Oton \cite{Ros-Oton} on the existence, regularity and qualitative properties of solutions to a class of elliptic nonlocal equations.\\
	
Motivated by the above works, it is natural to ask the following question:\\ 
	
\noindent Can we establish the regularity of viscosity solutions to \eqref{eq 0.1}?\\
	
The paper aims to answer this question affirmatively. To the best of our knowledge, we are not aware of such results for the mixed fully nonlinear local and nonlocal operators with variable-exponent gradient degeneracy.\\
The main challenges in this study are the following:
	\begin{enumerate}[label=\upshape{(C\arabic*)}, ref=(C\arabic*)]
		\item\label{C1} Due to the absence of the comparison principle, the classical Perron's method is not applicable here, which is quite standard and has been exploited in different contexts, see, for instance, \cite{Barles Busca, Filippis 1, Oza Heisenberg, Tyagi degenerate}.
		\item\label{C2} As mentioned earlier, the second challenge is that the non-local version of Lemma \ref{Lem Imbert} (see next) is not there. This prevents us to follow the classical approach of Imbert and Silvestre \cite{Imbert} in our case.
	\end{enumerate}
	
We follow the arguments of Birindelli and Demengel \cite{C1beta, Birin 1, Birin 4}, where the authors used the Ishii-Lions method \cite{Ishii} to establish the H\"older and Lipschitz regularity of the solutions to degenerate elliptic equations of type \eqref{a=1}. Further, using regularity results established in Theorem \ref{Main 1}, we get the compactness of the family of solutions to the approximated problem:
\begin{align}
		\mathlarger\varepsilon\big(\Delta u-(-\Delta)^{s}u\big)+|Du|^{\gamma(x)}\big(\a F\big(D^2u\big)-\b(-\Delta)^su\big)+b\cdot Du|Du|^{\gamma(x)}=f \text{ in }B_1,
	\end{align} 
	which further yields the existence of solution to \eqref{eq 0.1} by letting $\mathlarger\varepsilon\To 0^+.$ The similar arguments have been used effectively in \cite{Topp}. In this way, we overcome the challenge \ref{C1}. In dealing with \ref{C2} mentioned above, the order of the fractional Laplacian, $s\in({1}/{2},1)$ plays a crucial role. Here, we make use of the fact that $(-\Delta)^s$ approximates $-\Delta$ as $s\To 1^-,$ see \cite{Hitchhiker} for the details. \\
	
	
The novelty concerning the existing literature is due to the following:
	
	\begin{enumerate}
		\item [(i)] We deal with a variable-exponent degenerate mixed fully nonlinear local and nonlocal operator with first-order term $b\cdot Du.$
	\item [(ii)] The results are new even in the case when $\gamma(x)\equiv\gamma>0$ is a constant and $b\equiv 0$ in \eqref{eq 0.1}.
	\end{enumerate}
		
	Now, we list the assumptions on $F,$ $\gamma$ and $b,$ which we use in our main results:
	\begin{enumerate}[label=\upshape {(H\arabic*)}, ref=(H\arabic*)]
		\item\label{H3} (Ellipticity condition on $F$). $F\colon \mathcal{S}^N\To\mathbb{R}$ is uniformly ($\lambda,\,\Lambda$)-elliptic. In other words, for some parameters $\Lambda\geq \lambda>0,$ we have
		\begin{align}
			\lambda\|B\|\leq F(A+B)-F(A)\leq \Lambda\|B\|, 
		\end{align}
	for each $A,\,B\in \mathcal{S}^N$ with $B\geq 0.$ Also, $F(tA)=tF(A),$ $\forall t\geq 0$ and $A\in\mathcal{S}^N.$
		\item\label{H1} (Lower bound on $\gamma$). $\gamma\in C(\overline{B}_1,\mathbb{R})$ with ${\inf_{B_1}\gamma}\geq 0.$
		\item\label{H2} (Condition on $b$). $b:B_1\To\mathbb{R}^N$ is bounded Lipschitz with
		\begin{align}\label{b lip}
			(b(x)-b(y))\cdot (x-y)\leq 0, \text{ for }x,y\in B_1.		
		\end{align}
	\end{enumerate}
	
	The main results of this paper are the following theorems:
	
	\begin{thm}\label{Main 1}
		Let \ref{H3}, \ref{H1} and \ref{H2} hold. Let $f\in L^\infty(B_1,\mathbb{R})$ with $\|f\|_{\infty,B_1}\leq \varepsilon$ for some $\varepsilon>0$ and $s\in \big({1}/{2},1\big).$
		Then for each $p\in\mathbb{R}^N,$ any bounded viscosity solution $u$ of 
		\begin{align}
			\big|Du+p\big|^{\gamma(x)}\big(\a F(D^2u)-\b(-\Delta)^su\big)+b\cdot Du|Du+p|^{\gamma(x)}= f \text{ in }B_1,
		\end{align}
		with $\|u\|_{\infty,B_1}\leq 1$ is Lipschitz continuous with Lipschitz constant independent of $p.$ Moreover, the constant is uniformly bounded as $s\To 1^-.$ 
	\end{thm}
	
	\begin{thm}\label{main 2}
		Let \ref{H3}, \ref{H1} and \ref{H2} hold. Let $f\in L^\infty(B_1,\mathbb{R}).$ Then there exists $s_0\in({1}/{2},1)$ close enough to $1$ such that for any $s\in (s_0,1),$ every bounded viscosity solution of \eqref{eq 0.1} is $C^{1,\delta}$ for some $\delta\in (0,1)$ with
		\begin{align}
			[u]_{C^{1,\delta}(B_\frac{1}{2})}\leq	\widetilde{C}\left(\|u\|_{\infty,B_1}+\left(\|f\|_{\infty,B_1}^{\frac{2s-1}{1+\inf_{B_1}\gamma}}\right)\right)
		\end{align} 
		for some positive constant $\widetilde{C}.$
	\end{thm}
	
	We recall the following results (Corollary \ref{C1bt} to \ref{cor to}) due to several authors, which are closely related to Theorem \ref{main 2}.
	
	\begin{cor}[Remark 3.5 \cite{C1beta}]\label{C1bt}
		Let $\Omega\subset\mathbb{R}^N$ be a $C^2$ domain. Let $F$ be a uniformly elliptic operator and let $\gamma\geq 0$ be a constant. Let $h\in C(\overline{\Omega},\mathbb{R}^N)$ and $f\in C(\overline{\Omega},\mathbb{R}).$ Let $u$ be a viscosity solution of
		\begin{align}
			|Du|^\gamma(F(D^2u)+h\cdot Du)=f\text{ in }B_1	
		\end{align}
		with ${\text{osc}}_{B_1}u\leq 1$ and $\|f\|_{\infty,B_1}\leq \varepsilon,$ for some $\varepsilon\in [0,1].$ Then there exists $\tau,\delta\in(0,1)$ such that for all $n\in\mathbb{N},$ there exists $p_n\in\mathbb{R}^N$ such that
		\begin{align}
			\underset{B_{\tau^n}}{\text{osc}}~(u(x)-p_n\cdot x)\leq \tau^{n(1+\delta)}.
		\end{align} 
	\end{cor}
	
	\begin{cor}[Theorem 2.1 \cite{Bronzi}]
		Let $F$ be a uniformly elliptic operator. Let \ref{H1} hold and $f$ be bounded. Let $u$ be a viscosity solution of 
		\begin{align}
			|Du|^{\gamma(x)}F(D^2u)=f \text{ in }B_1.
		\end{align}
		Then $u$ is in $C^{1,\delta}$ for some $\delta\in(0,1).$	
	\end{cor}
		
		The following result is due to dos Prazeres and Topp \cite{Topp}. They considered a nonlocal nonlinear operator $\mathcal{I},$ defined using a class of nonlocal operators whose kernels are comparable to that of fractional Laplacian, i.e., $(-\Delta)^s.$ For the details, see \cite{Topp}.
		
		\begin{cor}[Theorem 1.2 \cite{Topp}]\label{cor to}
			Let $\gamma\geq 0$ be a constant and $f\in L^\infty(B_1,\mathbb{R})$. There exists $s_0\in\big({1}/{2},1\big)$ sufficiently close to $1$ such that for any $s\in (s_0,1),$ every bounded viscosity solution $u$ to
			\begin{align}
				|Du|^\gamma \mathcal{I}(u)=f \text{ in }B_1
			\end{align}
			is in $C^{1,\delta}$ for some $\delta\in (0,1).$
		\end{cor}
		
		\begin{rem}
			Let us recall from \cite{Imbert} that even in the case when $\lambda=\Lambda=1,$ i.e., $F(D^2u)=\Delta u,$ $\b=0,$ $b\equiv0$ and $\gamma(x)\equiv\gamma>0$, $f$ are constants, the best regularity one can expect for the solution of \eqref{eq 0.1} is $C^{1,\delta}.$
		\end{rem}
		
		\begin{rem}
			In the spirit of \cite{Arau}, in our case we have the following restriction on $\delta$$\colon$
			\begin{align}
				\delta< \displaystyle{\min}\left\{\overline{\delta},\frac{1}{1+\|\gamma^+\|_{\infty,B_1}}\right\},
			\end{align}
			where $\gamma^+$ denotes the positive part of $\gamma(x),$ and  $\overline{\delta}>0$ is such that the solutions of 
			\begin{align}
				\alpha F\big(D^2u\big)+\beta\Delta u+b\cdot Du=0 \text{ in }B_1
			\end{align}
			are $C^{1,\overline{\delta}}$ regular.
		\end{rem}
		The organization of this paper is as follows. In Section 2, we list the basic definitions and several important results from the literature. Section 3 is dedicated to the proof of our main results.

		\section{Preliminaries}
		Let $\mathcal{S}^N$ be the set of all real symmetric $N\times N$ matrices. Here, we have the usual partial ordering: $A\leq B$ in $\mathcal{S}^N$ means that $\langle A\xi,\xi\rangle \leq \langle B\xi, \xi\rangle$ for any $\xi \in \mathbb{R}^N$. In other words, $B-A$ is positive semidefinite.
		For $S\in \mathcal{S}^N$ and given two parameters $\Lambda\geq\lambda>0,$ extremal Pucci's operators are defined as follows:
		\begin{align}\label{M+}
			\mathcal{M}_{\lambda,\Lambda}^+(S):=\Lambda\displaystyle{\sum_{e_i\geq 0}e_i}+\lambda\displaystyle{\sum_{e_i< 0}e_i
			},
		\end{align}
	and
	\begin{align}\label{M-}
		\mathcal{M}_{\lambda,\Lambda}^-(S):=\Lambda\displaystyle{\sum_{e_i\leq 0}e_i}+\lambda\displaystyle{\sum_{e_i> 0}e_i
		},
	\end{align}
		where $\{e_i\}_{i=1}^N$ are the eigenvalues of $S.$ We mention that a $(\lambda,\Lambda)$-elliptic operator $F$ (see \ref{H3}) can also be defined in terms of Pucci's extremal operators as below.
		\begin{defn}
			An operator $F\colon\mathcal{S}^N\To \mathbb{R}$ is said to be uniformly elliptic with ellipticity constants $\lambda,\Lambda,$ i.e., $(\lambda,\Lambda)$-elliptic if the following holds$\colon$
			\begin{align}
				\mathcal{M}^-_{\lambda,\Lambda}(B)\leq F(A+B)-F(A)\leq \mathcal{M}^+_{\lambda,\Lambda}(B),
			\end{align}
		for each $A,B\in\mathcal{S}^N,$ with $B\geq 0.$
		\end{defn}
\begin{rem}\label{rem F}
	In particular, we have that
\begin{align}
	\mathcal{M}_{\lambda,\Lambda}^-(A)\leq F(A)\leq \mathcal{M}_{\lambda,\Lambda}^+(A),\,\forall A\in \mathcal{S}^N.
\end{align}
\end{rem}

		We now define the notion of viscosity solution, which we consider for \eqref{eq 0.1}. For $p\in\mathbb{R}^N,$ consider the auxiliary problem:
		\begin{align}\label{aux eq}
			\big|p+Du\big|^{\gamma(x)}\big(\a F\big(D^2u\big)-\b(-\Delta)^su\big)+b\cdot Du\big|p+Du\big|^{\gamma(x)}= f.
		\end{align}
		
		\begin{defn}\label{subsol}
			Let $u:\mathbb{R}^N\longrightarrow\mathbb{R}$ be an upper semicontinuous (USC) function in $\overline{B}_1.$ Then $u$ is called a \textit{viscosity subsolution} of \eqref{aux eq} if for any $x\in B_1$ and $C^2$ function $\varphi:\overline{B}_\rho\longrightarrow\mathbb{R},$ for some $\rho\in(0,1]$ such that $\varphi(x)=u(x),$ $\varphi(y)>u(y)$ for $y\in \overline{B}_\rho\setminus \{x\}$ and $D\varphi(x)\neq-p,$ we have
			\begin{align}
				\big|p+Dv(x)\big|^{\gamma(x)}\big(\a F\big(D^2v(x)\big)-\b(-\Delta)^sv(x)\big)+b(x)\cdot Dv(x)\big|p+Dv(x)\big|^{\gamma(x)}\geq f(x),
			\end{align}
			where
			\begin{align}
				v:=\begin{cases}
					\varphi &\text{ in }B_\rho,\\[1mm]
					u &\text{ in }\mathbb{R}^N\setminus B_\rho.
				\end{cases}
			\end{align}
			Moreover, we say $u$ satisfies $\big|p+Du\big|^{\gamma(x)}\big(\alpha F\big(D^2u\big)-\b(-\Delta)^su\big)+b\cdot Du\big|p+Du\big|^{\gamma(x)}\geq f \text{ in }B_1$ in the viscosity sense. 
		\end{defn}
		
		\begin{defn}
			Let $u:\mathbb{R}^N\longrightarrow\mathbb{R}$ be a lower semicontinuous (LSC) function in $\overline{B}_1.$ Then $u$ is called a \textit{viscosity supersolution} of \eqref{aux eq} if for any $x\in B_1$ and $C^2$ function $\psi:\overline{B}_\rho\longrightarrow\mathbb{R},$ for some $\rho\in(0,1]$ such that $\psi(x)=u(x),$ $\psi(y)<u(y)$ for $y\in \overline{B}_\rho\setminus \{x\}$ and $D\psi(x)\neq -p,$ we have
			\begin{align}
				\big|p+Dw(x)\big|^{\gamma(x)}\big(\a F\big(D^2w(x)\big)-\b(-\Delta)^sw(x)\big)+b(x)\cdot Dw(x)|p+Dw(x)|^{\gamma(x)}\leq f(x),
			\end{align}
			where
			\begin{align}
				w:=\begin{cases}
					\psi &\text{ in }B_\rho,\\[1mm]
					u &\text{ in }\mathbb{R}^N\setminus B_\rho.
				\end{cases}
			\end{align}
			Moreover, we say $u$ satisfies $\big|p+Du\big|^{\gamma(x)}\big(\a F\big(D^2u\big)-\b(-\Delta)^su\big)+b\cdot Du\big|p+Du\big|^{\gamma(x)}\leq f \text{ in }B_1$ in the viscosity sense.
		\end{defn}	
		
		We next	define the super/sub jets to demonstrate an equivalent definition of viscosity solution. Let us start by recalling the definitions of superjets and subjets.
		\begin{defn}\cite{Katzo}
			The \textit{superjet} $J^{2,+}u(x)$ of an USC function $u$ at $x\in \mathbb{R}^N$ consists of all the couples of the form $(p,X)\in \mathbb{R}^N\times \mathcal{S}^N$ for which
			\begin{align}
				u(x+z)\leq u(x)+p\cdot z+\frac{1}{2}\left\langle Xz,z\right\rangle+o(|z|^2). 
			\end{align} 
			Similarly, 
			the \textit{subjet} $J^{2,-}u(x)$ of a LSC function $u$ at $x\in \mathbb{R}^N$ consists of all the couples of the form $(p,X)\in\mathbb{R}^N\times \mathcal{S}^N$ for which
			\begin{align}
				u(x+z)\geq u(x)+p\cdot z+\frac{1}{2}\left\langle Xz,z\right\rangle+o(|z|^2). 
			\end{align}
		\end{defn}
		We next recall the definitions of the closures of the set-valued mappings, which will be required in the ensuing section.
		\begin{defn}\cite{User}
			The \textit{limiting superjet} $\overline{J}^{2,+}u(x)$ of a USC function $u$ is the set of couples of the form $(p,X)\in\mathbb{R}^N\times\mathcal{S}^N$ for which
			\begin{align}
				\exists (x_n,p_n,X_n)\longrightarrow(x,p,X) \text{ such that } (p_n,X_n) \text{ is a superjet of }u \text{ at }x_n \text{ and }u(x_n)\longrightarrow u(x). 
			\end{align}
			Similarly, the \textit{limiting subjet} $\overline{J}^{2,-}u(x)$ of an LSC function $u$ is the set of couples of the form $(p,X)\in\mathbb{R}^N\times\mathcal{S}^N$ for which
			\begin{align}
				\exists (x_n,p_n,X_n)\longrightarrow(x,p,X) \text{ such that } (p_n,X_n) \text{ is a subjet of }u \text{ at }x_n \text{ and }u(x_n)\longrightarrow u(x). 
			\end{align}
		\end{defn}
	 Imbert and Silvestre \cite{Imbert} proved the following result using an improvement of flatness technique:
		\begin{thm}[Theorem 1 \cite{Imbert}]\label{Thm Imbert}
			Let $\gamma\geq 0$ be a constant, $F$ be a uniformly elliptic operator with $F(0)=0$ and $f$ be bounded in $B_1.$ Then for any $U\Subset B_1,$ viscosity solution $u$ of
			\begin{align}
				|Du|^\gamma F(D^2u)=f \text{ in }B_1
			\end{align}   
			is in $C^{1,\delta}(U)$ with
			\begin{align}
				[u]_{C^{1,\delta}(U)}\leq C\left(\|u\|_{\infty,B_1}+\|f\|_{\infty,B_1}^{\frac{1}{1+\gamma}}\right).
			\end{align}
		\end{thm}
		The proof of the above theorem depends on the equicontinuity properties for the auxiliary equation:
		\begin{align}
			|p+Du|^\gamma F\big(D^2u\big)=f \text{ in } B_1,
		\end{align}
		in terms of $p\in\mathbb{R}^N.$ In the proof of Theorem \ref{Thm Imbert}, the following lemma plays a crucial role. 
		\begin{lem}[Lemma 6 \cite{Imbert}]\label{Lem Imbert}
			Let $u$ be a viscosity solution of 
			\begin{align}
				|p+Du|^\gamma F\big(D^2u\big)=0 \text{ in }B_1.
			\end{align}
			Then $u$ is a viscosity solution of 
			\begin{align}\label{cab}
				F\big(D^2u\big)=0 \text{ in }B_1.
			\end{align}
		\end{lem}
		Thanks to Caffarelli and Cabr\'e \cite{Cabre}, we have interior $C^{1,\delta}$ estimate for the solutions of \eqref{cab}. Imbert and Silvestre \cite{Imbert} also constructed a counterexample to justify that one can not expect better regularity than $C^{1,\delta}$ (see, Example 1 therein).
		
		We mention here that the presence of the nonlocal term in \eqref{eq 0.1} prevents us to get a nonlocal version of Lemma \ref{Lem Imbert}. In case, when $\a=0$ and $\b=1$, in \eqref{eq 0.1}, one may see the work of dos Prazeres and Topp \cite{Topp}, where they constructed a function, which immediately gives that for fractional Laplacian operator, $(-\Delta)^s,$ in place of $F,$ the analogue of Lemma \ref{Lem Imbert} does not hold, in general. For the sake of clarity, we mention below the same counterexample, which also works in our case, i.e., \eqref{eq 0.1} and is as follows. Consider a function $u:\mathbb{R}\longrightarrow\mathbb{R}$ given by  
		\begin{align}
			u(x)\coloneqq
			\begin{cases}
				x+1, &\text{ if }x\leq -1,\\[1mm]
				0, &\text{ if }x\in (-1,1),\\[1mm]
				x-1, &\text{ if }x\geq 1.
			\end{cases}
		\end{align}
		Let $s\in \big({1}/{2},1\big).$ Since $u\equiv 0$ in $(-1,1),$ so we have $$u_x=Du=0 \text{ and } u_{xx}=D^2u=0 \text{ in } (-1,1).$$
		
		\noindent One may see that
		\begin{align}
			|Du|^{\gamma(x)}\big(\a F\big(D^2u\big)-\b(-\Delta)^su\big)+b\cdot Du|Du|^{\gamma(x)}=0 &\text{ in } (-1,1).
		\end{align}
		Now, let $x\in(-1+\delta,1-\delta)\Subset(-1,1),$ for some $\delta>0.$  It yields
		\begin{align}
			\big(\a F\big(D^2u(x)\big)-\b(-\Delta)^su(x)\big)+b(x)\cdot Du(x)&=-\beta(-\Delta)^su(x)\\
			&=-\beta C(N,s) \,\text{ P.V.}\int_{\mathbb{R}^N}\frac{u(x)-u(y)}{|x-y|^{N+2s}}dy,
		\end{align}
		where $P.V.$ stands for the principal value and $C(N,s)$ is a constant depending on $N$ and $s,$ for the details see \cite{Hitchhiker}. It further gives
		\begin{align}\label{counter 1}\\
			\a F\big(D^2u(x)\big)-\b(-\Delta)^su(x)+b(x)\cdot Du(x)&=-\beta\,C(N,s) \,\displaystyle{\lim_{\varepsilon\to 0^+}}\left[\int_{-\infty}^{x-\varepsilon}\frac{u(x)-u(y)}{|x-y|^{1+2s}}dy+\int_{x+\varepsilon}^{\infty}\frac{u(x)-u(y)}{|x-y|^{1+2s}}dy\right]\\
			&=-\beta\,C(N,s) \,\displaystyle{\lim_{\varepsilon\to 0^+}}\left[-\int_{-\infty}^{x-\varepsilon}\frac{u(y)}{|x-y|^{1+2s}}dy-\int_{x+\varepsilon}^{\infty}\frac{u(y)}{|x-y|^{1+2s}}dy\right]\\
			&=-\beta\,C(N,s) \,\displaystyle{\lim_{\varepsilon\to 0^+}}\bigg[-\int_{-\infty}^{-1}\frac{u(y)}{|x-y|^{1+2s}}dy-\int_{-1}^{x-\varepsilon}\frac{u(y)}{|x-y|^{1+2s}}dy\\
			&\qquad\qquad\qquad\qquad\qquad-\int_{x+\varepsilon}^{1}\frac{u(y)}{|x-y|^{1+2s}}dy-\int_{1}^{\infty}\frac{u(y)}{|x-y|^{1+2s}}dy\bigg]\\
			&=-\beta\,C(N,s) \,\displaystyle{\lim_{\varepsilon\to 0^+}}\bigg[-\int_{-\infty}^{-1}\frac{y+1}{|x-y|^{1+2s}}dy-\int_{1}^{\infty}\frac{y-1}{|x-y|^{1+2s}}dy\bigg].
		\end{align}
		For the first term in the last step of R.H.S. of \eqref{counter 1}, $x-y\in(0,\infty)$ and for the second term $x-y\in (-\infty,0)$. This
		infers 
		\begin{align}\label{counter 2}\\
			\a F\big(D^2u(x)\big)-\b(-\Delta)^su(x)+b(x)\cdot Du(x)&=-\beta\,C(N,s)\,\displaystyle{\lim_{\varepsilon\rightarrow 0^+}}\bigg[-\int_{-\infty}^{-1}\frac{y+1}{(x-y)^{1+2s}}dy-\int_{1}^{\infty}\frac{y-1}{(y-x)^{1+2s}}dy\bigg]\\
			&=-\b C(N,s)\bigg[-\frac{2sy-x+2s-1}{(4s^2-2s)(x-y)^{2s}}\bigg|_{-\infty}^{-1}+\frac{2sy-x-2s+1}{(4s^2-2s)(y-x)^{2s}}\bigg|_{1}^{\infty}\bigg].
		\end{align} 
		Now, since $s\in\big({1}/{2},1\big),$ it gives that R.H.S. of \eqref{counter 2} is finite for $x\in (-1,1).$ Whereas, it blows up for $x\longrightarrow 1$ as well as for $x\longrightarrow -1.$ This counterexample elaborates on the effect of nonlocal term present with the gradient degeneracy. It is well known that $(-\Delta)^s\To -\Delta$ as $s\To 1^-,$ see \cite{Hitchhiker}. This property plays a crucial role in our analysis.
		
		Next, we state the non-local version of Ishii-Lions lemma \cite{Barles}, which we use in Theorem \ref{Main 1} to establish the Lipschitz regularity of viscosity solution to \eqref{eq 0.1}. It has also been exploited, for instance, in \cite{Barles 2, Barles 1}.
		\begin{lem}[Corollary 4.3.1 \cite{Barles}]\label{Ishii nonlocal}
			Let $F: \mathbb{R}^N\times \mathbb{R}^N\times\mathcal{S}^N\times\mathbb{R}\To\mathbb{R}$ be a continuous function satisfying the local and nonlocal degenerate ellipticity conditions. Let $u$ be a viscosity solution of
			\begin{align}
				F(x,Du,D^2u, \mathcal{I}[x,u])=0\text{ in }B_1.
			\end{align}
			Consider a function $v:B_1\times B_1\To\mathbb{R}$ given by
			\begin{align}
				v(x,y):=u(x)-u(y).
			\end{align} 
			Let $\Gamma\in C^2(B_1\times B_1)$ and $(x_1,y_1)\in B_1\times B_1$ be a point of global maxima of $v(x,y)-\Gamma(x,y)$ in $B_1\times B_1.$ Then, there exists $i>0$ small enough and matrices $X,\,Y\in\mathcal{S}^N$ such that
			\begin{align}
				F\left(x_1,D_x\Gamma(.,y_1),X,\mathcal{I}^{1,\delta}[x_1,\Gamma(.,y_1)],\mathcal{I}^{2,\delta}[x_1,D_x\Gamma(.,y_1),u]\right)\geq 0
			\end{align}
			and
			\begin{align}
				F\left(y_1,-D_y\Gamma(x_1,.),Y,\mathcal{I}^{1,\delta}[y_1,-\Gamma(x_1,.)],\mathcal{I}^{2,\delta}[y_1,-D_y\Gamma(x_1,.),u]\right)\leq 0,
			\end{align}
			where
			\begin{align}
				\mathcal{I}^{1,\delta}[x,\phi]=\int_{B_\delta}\big(\phi(x+y)-\phi(x)-D\phi(x)\cdot z\big)\mu(dy),
			\end{align}
			\begin{align}
				\mathcal{I}^{2,\delta}[x,p,w]=\int_{\mathbb{R}^N\setminus{B_\delta}}\big(w(x+y)-w(x)-\mathbbm{1}_{\{z\leq 1\}}z\cdot p\big)\mu(dy),
			\end{align}
			for $\delta\ll 1$ and suitable measure $\mu.$ Also, we have the following matrix inequality:
			\begin{align}
				\begin{bmatrix}
					X &0\\[1mm]
					0 &-Y
				\end{bmatrix}\leq \begin{bmatrix}
					D_{xx}^2\Gamma(x_1,y_1) &-D_{xy}^2\Gamma(x_1,y_1)\\[1mm]
					-D_{xy}^2\Gamma(x_1,y_1) &D_{yy}^2\Gamma(x_1,y_1)
				\end{bmatrix}+iI_{2N},	
			\end{align}
			where $I_{2N}$ denotes $2N\times 2N$ identity matrix.
		\end{lem}
		
		\section{Proofs of main results}
		\noindent \textbf{Proof of Theorem \ref{Main 1}.}
		We first state and prove the following two lemmas. The proof of Theorem \ref{Main 1} follows by combining these lemmas. 

		\begin{lem}\label{lemma 4}
			Let \ref{H3}, \ref{H1} and \ref{H2} hold and $s\in\big({1}/{2},1\big).$ Let $u$ be a viscosity solution of 
			\begin{align}\label{eq lem 4}
				\big|Du+p\big|^{\gamma(x)}\big(\a F(D^2u)-\b(-\Delta)^su\big)+b\cdot Du|Du+p|^{\gamma(x)}= f \text{ in }B_1
			\end{align}
			with ${\text{osc}}_{B_1}u\leq 1$ and $\|f\|_{\infty,B_1}\leq \varepsilon$ for some $\varepsilon>0.$ Then, there exists some $a>0$ such that if $|p|\geq {1}/{a},$ then $u$ is Lipschitz continuous with
			\begin{align}
				[u]_{C^{0,1}(B_{r})}\leq C=C(\lambda,\Lambda,N,r,\|f\|_{\infty,B_1},\|u\|_{\infty,B_1})
			\end{align}
			for each $r\in(0,{1}/{2}).$
		\end{lem}
		\noindent \textbf{Proof of Lemma \ref{lemma 4}.}
		Let $a_0={1}/{|p|}\in [0,a].$
		We use the viscosity method introduced by  Ishii and Lions \cite{Ishii}. Fix $r\in \big(0,1\big).$ For $x_0\in B_{{r}/{2}}$, define the auxiliary function
		\begin{align}\label{eq Phi}
			\Phi(x,y)&:=u(x)-u(y)-\Gamma(x,y)
		\end{align}
		with
		\begin{align}
			\Gamma(x,y)=L_1\phi(|x-y|)+L_2|x-x_0|^2+L_2|y-x_0|^2,
		\end{align}
		for $L_1>0,$ $L_2>0$ (to be chosen later) and
		\begin{align}\label{phi defn}
			\phi(t)=
			\begin{cases}
				t-\phi_0t^{1+\theta} &\text{ if }0<t\leq t_0\\[1mm]
				\phi(t_0) &\text{ if }t> t_0,
			\end{cases}
		\end{align}
		for the choice of $\theta\in(0,1)$ small enough. Here, we choose $\phi_0\in(0,{1}/{1+\theta})$ such that for some $t_0\geq 1$, $\phi(t_0)>0$ and $\phi_0(1+\theta)t^\theta<1$ for all $t\in (0,2).$ We look for $L_1$ and $L_2$ such that $$M:=\displaystyle{\sup_{B_1\times B_1}}~\Phi\leq 0.$$ Consider 
		\begin{align}
			L_2=\frac{16}{r^2},
		\end{align}
		for some fixed $r\in (0,1).$ We assume the contrary that $\Phi$ attains a positive maximum in $B_1\times B_1,$ i.e., $M>0$ for some $x_0\in B_{{r}/{2}}.$ This would further lead to a contradiction for $L_1$ large enough. Let $\Phi$ attains its maximum at some point $(x_1,y_1) \in \overline{B}_1\times\overline{B}_1.$ Using the continuity of $u$ and the fact that ${\text{osc}}_{B_1}\,u\leq 1,$ we have
		\begin{align}
			L_1\phi(|x_1-y_1|)+L_2|x_1-x_0|^2+L_2|y_1-x_0|^2\leq \underset{B_1}{\text{osc}}\,u\leq 1.
		\end{align}
		Since $\phi$ given by \eqref{phi defn} is positive, it infers
		\begin{align}
			L_2|x_1-x_0|^2+L_2|y_1-x_0|^2\leq 1.
		\end{align}
		This implies
		\begin{align}
			|x_1-x_0|^2+|y_1-x_0|^2\leq \frac{1}{L_2}=\frac{r^2}{16}. 
		\end{align}
		Using $x_0\in B_{{r}/{2}},$ it further concludes that $x_1, y_1$ belong to the interior of $B_1,$ in particular, in $B_{{3}/{4}}.$ Moreover, $x_1\neq y_1,$ otherwise $M\leq 0.$
		Also, it is easy to see that
		\begin{align}
			\phi'(t)=\begin{cases}
				1-(1+\theta) t^{\theta}\phi_0 &\text{ if }0< t\leq t_0,\\[1mm]
				0 &\text{ if }t>t_0.
			\end{cases}
		\end{align}
		Observe that $x\mapsto u(x)-\Gamma(x,y_1)$ attains a maximum in $B_1$ at $x_1.$ Similarly, $y\mapsto -u(y)-\Gamma(x_1,y)$ attains a maximum in $B_1$ at $y_1.$ In other words, $u(y)-\big(-\Gamma(x_1,y)\big)$ attains a minimum in $B_1$ at $y_1.$ By the definition of viscosity solution, $\Gamma(x_1,.)$ and $-\Gamma(.,y_1)$ work as test functions for viscosity subsolution and supersolution, respectively. 
		Also, we have
		\begin{align}
			D_x\Gamma(x_1,y_1)&=L_1\phi'\big(|x_1-y_1|\big)\frac{(x_1-y_1)}{|x_1-y_1|}+2L_2(x_1-x_0),
		\end{align}
		and
		\begin{align}
			D_y\Gamma(x_1,y_1)&=-L_1\phi'(|x_1-y_1|)\frac{(x_1-y_1)}{|x_1-y_1|}+2L_2(y_1-x_0).
		\end{align}
		From the above equations, note that
		\begin{align}\label{est D Gamma}
			\big|D_x\Gamma(x_1-y_1)-\big(-D_y\Gamma(x_1,y_1)\big)\big|\leq 4L_2.	
		\end{align}
		Also,
		\begin{align}\label{eq cL1}
			|D_x\Gamma(x_1,y_1)|\leq L_1|\phi'(x_1-y_1)|+2L_2\leq cL_1,	
		\end{align}
		for some universal constant $c>0.$ Next, our aim is to write two viscosity inequalities, for which we construct limiting subjet $(D_x\Gamma(.,y_1),X)$ of $u$ at $x_1$ and a limiting superjet $(-D_y\Gamma(x_1,.),Y)$ of $u$ at $y_1$ using nonlocal version of Jensen-Ishii's Lemma (Lemma \ref{Ishii nonlocal}). We have the following matrix inequality for small enough $i>0$$\colon$
		\begin{align}
			\begin{bmatrix}
				X &0\\[1mm]
				0 &-Y
			\end{bmatrix}
			&\leq \begin{bmatrix}
				L_1D^2\left(\phi(|.|)\right)(x_1-y_1)+2L_2I_N &-L_1D^2\left(\phi(|.|)\right)(x_1-y_1)\\[1mm]
				-L_1D^2\left(\phi(|.|)\right)(x_1-y_1) &L_1D^2\left(\phi(|.|)\right)(x_1-y_1)+2L_2I_N
			\end{bmatrix}+iI_{2N}\\
			&= \begin{bmatrix}
				L_1D^2\left(\phi(|.|)\right)(x_1-y_1) &-L_1D^2\left(\phi(|.|)\right)(x_1-y_1)\\[1mm]
				-L_1D^2\left(\phi(|.|)\right)(x_1-y_1) &L_1D^2\left(\phi(|.|)\right)(x_1-y_1)
			\end{bmatrix}+\begin{bmatrix}
				2L_2I_N &0\\[1mm]
				0 &2L_2I_N
			\end{bmatrix}+iI_{2N}\\
			&=\begin{bmatrix}
				Z &-Z\\[1mm]
				-Z &Z
			\end{bmatrix}+2L_2\begin{bmatrix}
				I_N &0\\[1mm]
				0 &I_N
			\end{bmatrix}+iI_{2N},
		\end{align}
		where \begin{align}
			Z=L_1\bigg(\phi''(|x_1-y_1|)\frac{(x_1-y_1)}{|x_1-y_1|}\otimes\frac{(x_1-y_1)}{|x_1-y_1|}+\frac{\phi'(|x_1-y_1|)}{|x_1-y_1|}\left(I_N-\frac{(x_1-y_1)}{|x_1-y_1|}\otimes\frac{(x_1-y_1)}{|x_1-y_1|}\right)\bigg),
		\end{align}
		and $I_N,$ $I_{2N}$ denote $N\times N$ and $2N\times 2N$ identity matrices, respectively. Next, applying the above inequality to any vector $(\xi,\xi)$ with $|\xi|=1,$ we have that
		\begin{align}
			\left\langle (X-Y)\xi,\xi\right\rangle\leq 4L_2+2i.
		\end{align}
		Also, applying the matrix inequality to the vector $\left(\frac{(x_1-y_1)}{|x_1-y_1|},-\frac{(x_1-y_1)}{|x_1-y_1|}\right),$ we get
		\begin{align}
			\left\langle (X-Y)\frac{(x_1-y_1)}{|x_1-y_1|},\frac{(x_1-y_1)}{|x_1-y_1|} \right\rangle&\leq 4\left\langle Z\frac{(x_1-y_1)}{|x_1-y_1|},\frac{(x_1-y_1)}{|x_1-y_1|}\right\rangle+(4L_2+2i)\left\langle I_N\frac{(x_1-y_1)}{|x_1-y_1|}, \frac{(x_1-y_1)}{|x_1-y_1|}\right\rangle\\
			&\leq \bigg(-4(1+\theta)\theta L_1|x_1-y_1|^{\theta-1}\phi_0+4L_2+2i\bigg)\left|\frac{x_1-y_1}{|x_1-y_1|}\right|^2\\
			&\leq \bigg(-4(1+\theta)\theta L_12^{\theta-1}\phi_0+4L_2+2i\bigg)\left|\frac{x_1-y_1}{|x_1-y_1|}\right|^2. 
		\end{align}
		It infers that the matrix $X-Y$ has atleast one eigenvalue less than $\left(-4(1+\theta)\theta L_12^{\theta-1}\phi_0+4L_2+i\right).$ This quantity will be negative for sufficiently large values of $L_1.$ Let $\{e_i\}_{i=1}^N$ be denote the eigenvalues of the matrix $X-Y.$ We observe that, by definition,
		\begin{align}\label{in eq 1}
			\mathcal{M}_{\lambda,\Lambda}^+(X-Y)&=\Lambda\sum_{e_i\geq 0}e_i+\lambda\sum_{e_i<0}e_i\\
			&\leq \lambda\left(-4(1+\theta)\theta L_12^{\theta-1}\phi_0+4L_2+2i\right)+\Lambda(N-1)(4L_2+2i).
		\end{align}
		We choose $a:={1}/{5cL_1},$ where $c$ is the same constant appearing in \eqref{eq cL1}. This permits
		\begin{align}\label{p++}
			|{p}+D_x\Gamma(x_1,y_1)|&\geq |p|-|D_x\Gamma(x_1,y_1)|\\
			&\geq 5cL_1-cL_1\\
			&=4cL_1\\
			&=\frac{4}{5a}.
		\end{align}
Similarly,
\begin{align}\label{p--}
	|p-D_y\Gamma(x_1,y_1)|&\geq|p|-D_y\Gamma(x_1,y_1)\\
	&\geq \frac{4}{5a}.
\end{align}
		Next, we write the viscosity inequalities for $u$ at $x_1$ and $y_1$$\colon$
		\begin{align}\label{in eq1}
			{f}(x_1)&\leq \big|{p}+D_x\Gamma(x_1,y_1)\big|^{\gamma(x)}\big(\a F(X)-\b(-\Delta)^s{u}(x_1)+b(x_1)\cdot D_x\Gamma(x_1,y_1)\big)\\&=
			\big|{p}+D_x\Gamma(x_1,y_1)\big|^{\gamma(x)}\bigg(\a F(X)+\b C(N,s)\int_{B_\delta}\frac{\big(\Gamma(x_1+z,y_1)-\Gamma(x_1,y_1)-{\mathbbm{1}_{\{|z|\leq 1\}}}z\cdot D_x\Gamma(x_1,y_1)\big)}{|z|^{N+2s}}dz\\
			&\quad\qquad\qquad\qquad\qquad\qquad+\b C(N,s)\int_{\mathbb{R}^N\setminus B_\delta}\frac{\big(u(x_1+z)-u(x_1)-{\mathbbm{1}_{\{z\leq 1\}}}z\cdot D_x\Gamma(x_1,y_1)\big)}{|z|^{N+2s}}dz+b(x_1)\cdot D_x\Gamma(x_1,y_1)\bigg),
		\end{align}
		and
		\begin{align}\label{in eq2}
			{f}(y_1)&\geq \big|{p}-D_y\Gamma(x_1,y_1)\big|^{\gamma(x)}\big(\a F(Y)-\b(-\Delta)^s{u}(y_1)-b(y_1)\cdot D_y\Gamma(x_1,y_1)\big)\\
			&=\big|{p}-D_y\Gamma(x_1,y_1)\big|^{\gamma(x)}\bigg(\a F(Y)+\b C(N,s)\int_{B_\delta}\frac{\big(-\Gamma(x_1,y_1+z)+\Gamma(x_1,y_1)+{\mathbbm{1}_{\{|z|\le 1\}}}z\cdot D_y\Gamma(x_1,y_1)\big)}{|z|^{N+2s}}dz\\
			&\quad\qquad\qquad\qquad\qquad\qquad+\b C(N,s)\int_{\mathbb{R}^N\setminus B_\delta}\frac{\big(u(y_1+z)-u(y_1)+{\mathbbm{1}_{\{|z|\leq 1\}}}z\cdot D_y\Gamma(x_1,y_1)\big)}{|z|^{N+2s}}dz-b(y_1)\cdot D_y\Gamma(x_1,y_1)\bigg)
		\end{align}
		for small enough $\delta>0.$ Above inequalities further infer
		\begin{align}\label{in eq11}
			\quad\quad\quad \big|{p}+D_x\Gamma(x_1,y_1)\big|^{\gamma(x)}&\bigg(\a  F(X)+\b C(N,s)\int_{B_\delta}\frac{\big(\Gamma(x_1+z,y_1)-\Gamma(x_1,y_1)-{\mathbbm{1}_{\{|z|\leq 1\}}}z\cdot D_x\Gamma(x_1,y_1)\big)}{|z|^{N+2s}}dz\\
			&\quad+\b C(N,s)\int_{\mathbb{R}^N\setminus B_\delta}\frac{\big(u(x_1+z)-u(x_1)-{\mathbbm{1}_{\{z\leq 1\}}}z\cdot D_x\Gamma(x_1,y_1)\big)}{|z|^{N+2s}}dz\bigg)\\
			&\quad+\big|{p}+D_x\Gamma(x_1,y_1)\big|^{\gamma(x)}b(x_1)\cdot D_x\Gamma(x_1,y_1)\geq -\|{f}\|_{\infty,B_1},
		\end{align}
		and
		\begin{align}\label{in eq21}
			\quad\quad\quad \big|{p}-D_y\Gamma(x_1,y_1)\big|^{\gamma(x)}&\bigg(\a F(Y)+\b C(N,s)\int_{B_\delta}\frac{\big(-\Gamma(x_1,y_1+z)+\Gamma(x_1,y_1)+{\mathbbm{1}_{\{|z|\leq 1\}}}z\cdot D_y\Gamma(x_1,y_1)\big)}{|z|^{N+2s}}dz\\
			&\quad+\b C(N,s)\int_{\mathbb{R}^N\setminus B_\delta}\frac{\big(u(y_1+z)-u(y_1)+{\mathbbm{1}_{\{|z|\le 1\}}}z\cdot D_y\Gamma(x_1,y_1)\big)}{|z|^{N+2s}}dz\bigg)\\
			&\quad-\big|{p}-D_y\Gamma(x_1,y_1)\big|^{\gamma(x)}b(y_1)\cdot D_y\Gamma(x_1,y_1)\leq \|{f}\|_{\infty,B_1}. 
		\end{align}
		Subtracting \eqref{in eq21} from \eqref{in eq11} yields
		\begin{align}\label{in eq 3}
			-2\|{f}\|_{\infty,B_1}&\leq\big|{p}+D_x\Gamma(x_1,y_1)\big|^{\gamma(x)}\a F(X)-\big|{p}-D_y\Gamma(x_1,y_1)\big|^{\gamma(x)}\a F(Y)\\
			&\qquad\qquad+\b C(N,s)\big|{p}+D_x\Gamma(x_1,y_1)\big|^{\gamma(x)}\bigg(\int_{B_\delta}\frac{\big(\Gamma(x_1+z,y_1)-\Gamma(x_1,y_1)-{\mathbbm{1}_{\{|z|\leq 1\}}}z\cdot D_x\Gamma
				(x_1,y_1)\big)}{|z|^{N+2s}}dz\\
			&\qquad\qquad\qquad\qquad\qquad\qquad\qquad\qquad+\int_{\mathbb{R}^N\setminus B_\delta}\frac{\big(u(x_1+z)-u(x_1)-{\mathbbm{1}_{\{z\leq 1\}}}z\cdot D_x\Gamma(x_1,y_1)\big)}{|z|^{N+2s}}dz\bigg)\\
			&\qquad\qquad-\b C(N,s)\big|{p}-D_y\Gamma(x_1,y_1)\big|^{\gamma(x)}\bigg(\int_{B_\delta}\frac{\big(-\Gamma(x_1,y_1+z)+\Gamma(x_1,y_1)+{\mathbbm{1}_{\{|z|\le 1\}}}z\cdot D_y\Gamma(x_1,y_1)\big)}{|z|^{N+2s}}dz\\
			&\qquad\qquad\qquad\qquad\qquad\qquad\qquad\qquad+\int_{\mathbb{R}^N\setminus B_\delta}\frac{\big(u(y_1+z)-u(y_1)+{\mathbbm{1}_{\{|z|\leq 1\}}}z\cdot D_y\Gamma(x_1,y_1)\big)}{|z|^{N+2s}}dz\bigg).	
		\end{align} 
		\noindent Here, we use the idea, which is quite standard and has been exploited in \cite{Barles 2, Barles 1, Topp}. For $|z|\leq{1}/{4},$ we have $x_1+z,\,y_1+z\in \overline{B}_1.$
		Define
		\begin{align}
			I_1&\coloneqq\b C(N,s)\int_{B_{\frac{1}{4}}\setminus B_\delta}\frac{\big(u(x_1+z)-u(x_1)-{\mathbbm{1}_{\{|z|\leq 1\}}}z\cdot D_x\Gamma(x_1,y_1)\big)}{|z|^{N+2s}}dz\\
			&\qquad\qquad\qquad-\b C(N,s)\int_{B_{\frac{1}{4}}\setminus B_\delta}\frac{\big(u(y_1+z)-u(y_1)+{\mathbbm{1}_{\{|z|\le 1\}}}z\cdot D_y\Gamma(x_1,y_1)\big)}{|z|^{N+2s}}dz\\
			I_2&\coloneqq\b C(N,s)\int_{B_\delta}\frac{\big(\Gamma(x_1+z,y_1)-\Gamma(x_1,y_1)-{\mathbbm{1}_{\{|z|\leq 1\}}}z\cdot D_x\Gamma(x_1,y_1)\big)}{|z|^{N+2s}}dz\\
			&\qquad\qquad\qquad-\b C(N,s)\int_{B_\delta}\frac{\big(-\Gamma(x_1,y_1+z)+\Gamma(x_1,y_1)+{\mathbbm{1}_{\{|z|\le 1\}}}z\cdot D_y\Gamma(x_1,y_1)\big)}{|z|^{N+2s}}dz.	
		\end{align}
		We define a set 
		\begin{align}
			\mathcal{C}=\bigg\{z\in B_\rho:(1-\eta)|z||x_1-y_1|\leq |(x_1-y_1)\cdot z|\bigg\}
		\end{align}
		for some fixed $\rho\in(0,{1}/{4}).$ Consider
		\begin{align}
			I_3:=&\b C(N,s)\int_{\mathcal{C}}\frac{\big(u(x_1+z)-u(x_1)-{\mathbbm{1}_{\{|z|\leq 1\}}}z\cdot D_x\Gamma(x_1,y_1)\big)}{|z|^{N+2s}}dz\\
			&-\b C(N,s)\int_{\mathcal{C}}\frac{\big(u(y_1+z)-u(y_1)+{\mathbbm{1}_{\{|z|\le 1\}}}z\cdot D_y\Gamma(x_1,y_1)\big)}{|z|^{N+2s}}dz.
		\end{align}
		Since $(x_1,y_1)$ is a point of maxima for $\Phi$ in $B_1\times B_1$ and $M$ is the maximum, so we have
		\begin{align}
			u(x_1+z)-u(y_1+\overline{z})-\Gamma(x_1+z,y_1+\overline{z})\leq u(x_1)-u(y_1)-\Gamma(x_1,y_1),
		\end{align}
		for any $z,\,\overline{z}\in B_{{1}/{4}}.$ 
		It further displays
		\begin{align}\label{zz'}
			u(x_1+z)-u(x_1)-D_x\Gamma(x_1,y_1)\cdot z&\leq u(y_1+\overline{z})-u(y_1) +D_y\Gamma(x_1,y_1)\cdot\overline{z}\\
			&\quad+\Gamma(x_1+z,y_1+\overline{z})-\Gamma(x_1,y_1)-D\Gamma(x_1,y_1)\cdot (z,\overline{z}),
		\end{align}
		where $$D\Gamma(x_1,y_1)=\big(D_x\Gamma(x_1,y_1),D_y\Gamma(x_1,y_1)\big).$$ Taking, in particular $\overline{z}=0$ yields
		\begin{align}
			u(x_1+z)-u(x_1)-D_x\Gamma(x_1,y_1)\cdot z\leq \Gamma(x_1+z,y_1)-\Gamma(x_1,y_1)-D_x\Gamma(x_1,y_1)\cdot z.
		\end{align}
		Similarly, putting $z=0$ infers
		\begin{align}
			-\big(u(y_1+\overline{z})-u(y_1) +D_y\Gamma(x_1,y_1)\cdot\overline{z}\big)&\leq \Gamma(x_1,y_1+\overline{z})-\Gamma(x_1,y_1)-D_y\Gamma(x_1,y_1)\cdot \overline{z}.
		\end{align}
		Using above inequalities, we have
		\begin{align}
			I_3&= \beta C(N,s)\int_{\mathcal{C}}\frac{\big(u(x_1+z)-u(x_1)-z\cdot D_x\Gamma(x_1,y_1)\big)}{|z|^{N+2s}}dz\\
			&\qquad-\beta C(N,s)\int_{\mathcal{C}}\frac{\big(u(y_1+z)-u(y_1)+z\cdot D_y\Gamma(x_1,y_1)\big)}{|z|^{N+2s}}dz\\
			&\leq \beta C(N,s)\int_{\mathcal{C}}\frac{\big(\Gamma(x_1+z,y_1)-\Gamma(x_1,y_1)-z\cdot D_x\Gamma(x_1,y_1)\big)}{|z|^{N+2s}}dz\\
			&\qquad+\beta C(N,s)\int_{\mathcal{C}}\frac{\big(\Gamma(x_1,y_1+z)-\Gamma(x_1,y_1)-z\cdot D_y\Gamma(x_1,y_1)\big)}{|z|^{N+2s}}dz.
		\end{align}
		Now, since $$\Gamma(x,y)=L_1\phi(|x-y|)+L_2|x-x_0|^2+L_2|y-x_0|^2.$$ It yields
		\begin{align}
			I_3
			&\leq L_1\beta C(N,s)\int_{\mathcal{C}}\frac{\big(\phi(|x_1-y_1+z|)-\phi(|x_1-y_1|)-z\cdot D\phi(|x_1-y_1|)\big)}{|z|^{N+2s}}dz\\
			&\quad+L_2\beta C(N,s)\int_{\mathcal{C}}\frac{\big(|x_1+z-x_0|^2-|x_1-x_0|^2-z\cdot 2(x-x_0)\big)}{|z|^{N+2s}}dz\\
			&\quad+L_1\beta C(N,s)\int_{\mathcal{C}}\frac{\big(\phi(|x_1-y_1-z|)-\phi(|x_1-y_1|)+z\cdot D\phi(|x_1-y_1|)\big)}{|z|^{N+2s}}dz\\
			&\quad+L_2\beta C(N,s)\int_{\mathcal{C}}\frac{\big(|y_1+z-x_0|^2-|y_1-x_0|^2-z\cdot 2(y-x_0)\big)}{|z|^{N+2s}}dz\\
			&=L_1\beta C(N,s)\int_{\mathcal{C}}\frac{\big(\phi(|x_1-y_1+z|)-\phi(|x_1-y_1|)-z\cdot D\phi(|x_1-y_1|)\big)}{|z|^{N+2s}}dz+L_2\beta C(N,s)\int_{\mathcal{C}}\frac{|z|^2}{|z|^{N+2s}}dz\\
			&\qquad +L_1\beta C(N,s)\int_{\mathcal{C}}\frac{\big(\phi(|x_1-y_1-z|)-\phi(|x_1-y_1|)+z\cdot D\phi(|x_1-y_1|)\big)}{|z|^{N+2s}}dz+L_2\beta C(N,s)\int_{\mathcal{C}}\frac{|z|^2}{|z|^{N+2s}}dz.
		\end{align}
		Further, using a second-order Taylor expansion in the above integrals, we get
		\begin{align}\label{I2 est}
			\\I_3&\leq \frac{L_1\beta C(N,s)}{2}\left(\int_{\mathcal{C}}\sup_{t\in(0,1)}\frac{D^2\phi(|x_1-y_1+tz|)z\cdot z}{|z|^{N+2s}}dz+\int_{\mathcal{C}}\sup_{t\in(-1,0)}\frac{D^2\phi(|x_1-y_1+tz|)z\cdot z}{|z|^{N+2s}}dz\right)\\
			&\qquad+2L_2\beta C(N,s)\int_{\mathcal{C}}\frac{|z|^2}{|z|^{N+2s}}dz\\	&\leq\frac{L_1\beta C(N,s)}{2}\int_{\mathcal{C}}\sup_{t\in(-1,1)}\frac{D^2\phi(|x_1-y_1+tz|)z\cdot z}{|z|^{N+2s}}dz+2L_2\beta C(N,s)\int_{\mathcal{C}}\frac{|z|^2}{|z|^{N+2s}}dz\\
			&=2L_2\beta C(N,s)\int_{\mathcal{C}}\frac{|z|^2}{|z|^{N+2s}}dz+\frac{L_1\beta C(N,s)}{2}\times\\
			&\qquad\int_{\mathcal{C}}\sup_{t\in(-1,1)}\frac{\bigg(\phi''(|x_1-y_1+tz|)\frac{(x_1-y_1+tz)}{|x_1-y_1+tz|}\otimes\frac{(x_1-y_1+tz)}{|x_1-y_1+tz|}+\frac{\phi'(|x_1-y_1+tz|)}{|x_1-y_1+tz|}\left(I-\frac{(x_1-y_1+tz)}{|x_1-y_1+tz|}\otimes\frac{(x_1-y_1+tz)}{|x_1-y_1+tz|}\right)\bigg)z\cdot z}{|z|^{N+2s}}dz\\
			&=2L_2\beta C(N,s)\int_{\mathcal{C}}\frac{|z|^2}{|z|^{N+2s}}dz+\frac{L_1\beta C(N,s)}{2}\times\\
			&\qquad\int_{\mathcal{C}}\sup_{t\in(-1,1)}\frac{\phi''(|x_1-y_1+tz|)\bigg|\left(\frac{x_1-y_1+tz}{|x_1-y_1+tz|}\right)\cdot z\bigg|^2+\frac{\phi'(|x_1-y_1+tz|)}{|x_1-y_1+tz|}\left(|z|^2-\bigg|\left(\frac{x_1-y_1+tz}{|x_1-y_1+tz|}\right)\cdot z\bigg|^2\right)}{|z|^{N+2s}}dz.\\
		\end{align}
		Let $$\rho=\nu_0|x_1-y_1|<|x_1-y_1| \text{ for } \nu_0\in(0,1),$$ small enough to be chosen later. Now, we estimate $|x_1-y_1+tz|$ and $(x_1-y_1+tz)\cdot z$ on the set $\mathcal{C}$ as follows:
		\begin{align}
			|x_1-y_1+tz|&\leq |x_1-y_1|+|t||z|\leq|x_1-y_1|+\rho=|x_1-y_1|(1+\nu_0).\\
			|x_1-y_1+tz|&\geq |x_1-y_1|-|t||z|\geq |x_1-y_1|-\rho=|x_1-y_1|(1-\nu_0).\\ 
			|(x_1-y_1+tz)\cdot z|&=|(x_1-y_1)\cdot z+t|z|^2|
			\\
			&\geq|(x_1-y_1)\cdot z|-\rho|z|\\
			&\geq(1-\eta)|z||x_1-y_1|-\rho|z|\\
			&=(1-\eta)|z||x_1-y_1|-|x_1-y_1|\nu_0|z|\\
			&=(1-\eta-\nu_0)|x_1-y_1||z|.
		\end{align}
		Using the above inequalities, we get 
		\begin{align}
			\left|\frac{x_1-y_1+tz}{|x_1-y_1+tz|}\cdot z\right|&\geq \frac{1}{|x_1-y_1+tz|}(1-\eta-\nu_0)|x_1-y_1||z|\\
			&\geq \frac{(1-\eta-\nu_0)|x_1-y_1||z|}{|x_1-y_1|(1+\nu_0)}\\
			&=\frac{(1-\eta-\nu_0)|z|}{(1+\nu_0)},
		\end{align}
		for $t\in(-1,1).$ Using \eqref{phi defn} in \eqref{I2 est} infers
		\begin{align}\label{est2 I2}
			I_3&\leq \frac{L_1\beta}{2}\int_{\mathcal{C}}\sup_{t\in(-1,1)}\frac{-(1+\theta)\theta|x_1-y_1+tz|^{\theta-1}\phi_0\left(\frac{(1-\eta-\nu_0)^2|z|^2}{(1+\nu_0)^2}\right)}{|z|^{N+2s}}dz\\
			&\qquad+\frac{L_1\beta}{2}\int_{\mathcal{C}}\sup_{t\in(-1,1)}\frac{\frac{1}{|x_1-y_1+tz|}\bigg({1}-(1+\theta)\phi_0|x_1-y_1+tz|^{\theta}\bigg)\left(|z|^2-\bigg|\left(\frac{x_1-y_1+tz}{|x_1-y_1+tz|}\right)\cdot z\bigg|^2\right)}{|z|^{N+2s}}dz\\
			&\qquad+2L_2\beta C(N,s)\int_{\mathcal{C}}\frac{|z|^2}{|z|^{N+2s}}dz\\
			&\leq L_1\beta C(N,s)\int_{\mathcal{C}}\sup_{t\in(-1,1)}\frac{-(1+\theta)\theta|x_1-y_1+tz|^{\theta-1}\phi_0\bigg(\frac{(1-\eta-\nu_0)^2|z|^2}{(1+\nu_0)^2}\bigg)}{|z|^{N+2s}}dz\\
			&\qquad+L_1\beta C(N,s)\int_{\mathcal{C}}\sup_{t\in(-1,1)}\frac{\frac{1}{|x_1-y_1+tz|}\bigg(1-(1+\theta)\phi_0|x_1-y_1+tz|^{\theta}\bigg)\left(1-\frac{(1-\eta-\nu_0)^2}{(1+\nu_0)^2}\right)|z|^2}{|z|^{N+2s}}dz\\
			&\qquad+2L_2\beta C(N,s)\int_{\mathcal{C}}\frac{|z|^2}{|z|^{N+2s}}dz.
		\end{align}
		Note that $$\frac{(1-\eta-\nu_0)^2}{(1+\nu_0)^2}\leq 1.$$ Taking $\eta=c_1|x_1-y_1|^{2\theta}$ and $\nu_0=c_1|x_1-y_1|^{\theta}$ for some constant $c_1>0,$ we obtain $\frac{1-\eta-\nu_0}{1+\nu_0}\simeq 1,$ see {\rm \cite[Remark 7]{Barles 2}}. Hence the term $I_3$ can be made as negative as we please. 
		Also, we have
		\begin{align}
			I_2=\beta C(N,s)\int_{B_{\frac{1}{4}}\setminus B_\delta}\frac{\big(\Gamma(x_1+z,y_1+z)-\Gamma(x_1,y_1)-D\Gamma(x_1,y_1)\cdot (z,z)\big)}{|z|^{N+2s}}dz.
		\end{align}
		We now use a second order Taylor expansion for $\Gamma$ and get
		\begin{align}
			I_2\leq \frac{\beta C(N,s)}{2} \int_{B_\frac{1}{4}\setminus B_{\delta}}\frac{D^2\Gamma(x_1,y_1)(z,z)\cdot (z,z)}{|z|^{N+2s}}dz.
		\end{align}
		Further, taking $L_1$ large enough, we get using the concavity of $\phi$ that $I_2\leq 0.$
		Similarly, we have $I_1\leq 0.$ Combining the above estimates, we have that the sum $I_1+I_2+I_3$ can be made as negative as we want, let us say $-P,$ where $P=P(L_1)$ is a large enough positive number. On the other hand, since $u$ is bounded and the kernel is uniformly bounded away from the origin, this yields
	    that the integral 
	    \begin{align}
	    	\int_{\mathbb{R}^N\setminus{B_{1/4}}}\frac{u(x_1+z)-u(x_1)}{|z|^{N+2s}}dz-\int_{\mathbb{R}^N\setminus{B_{1/4}}}\frac{u(y_1+z)-u(y_1)}{|z|^{N+2s}}dz
	    \end{align}
	     is uniformly bounded. Using the above estimates, we deduce
		\begin{align}
			\frac{{f}(x_1)}{\big|{p}+D_x\Gamma(x_1,y_1)\big|^{\gamma(x_1)}}&\leq \bigg(\a F(X)+\b C(N,s)\int_{B_\delta}\frac{\big(\Gamma(x_1+z,y_1)-\Gamma(x_1,y_1)-{\mathbbm{1}_{\{|z|\leq 1\}}}z\cdot D_x\Gamma(x_1,y_1)\big)}{|z|^{N+2s}}dz\\
			&\qquad+\b C(N,s)\int_{\mathbb{R}^N\setminus B_\delta}\frac{\big(u(x_1+z)-u(x_1)-{\mathbbm{1}_{\{z\leq 1\}}}z\cdot D_x\Gamma(x_1,y_1)\big)}{|z|^{N+2s}}dz\bigg)+b(x_1)\cdot D_x\Gamma(x_1,y_1)\\	
			&=\bigg(\a F(X)+\b C(N,s)\int_{B_\delta}\frac{\big(\Gamma(x_1+z,y_1)-\Gamma(x_1,y_1)-{\mathbbm{1}_{\{|z|\leq 1\}}}z\cdot D_x\Gamma(x_1,y_1)\big)}{|z|^{N+2s}}dz\\
			&\qquad+\b C(N,s)\int_{\mathbb{R}^N\setminus B_\delta}\frac{\big(u(x_1+z)-u(x_1)-{\mathbbm{1}_{\{z\leq 1\}}}z\cdot D_x\Gamma(x_1,y_1)\big)}{|z|^{N+2s}}dz\bigg)\\	
			&\qquad+\big(b(x_1)-b(y_1)\big)\cdot D_x\Gamma(x_1,y_1)-b(y
			_1)\cdot D_y\Gamma(x_1,y_1)+b(y_1)\cdot \big(D_x\Gamma(x_1,y_1)+D_y\Gamma(x_1,y_1)\big)
			\end{align}
		On using \eqref{est D Gamma}, we have		
			\begin{align}
			\frac{{f}(x_1)}{\big|{p}+D_x\Gamma(x_1,y_1)\big|^{\gamma(x_1)}}&\leq \bigg(\a F(X)+\b C(N,s)\int_{B_\delta}\frac{\big(\Gamma(x_1+z,y_1)-\Gamma(x_1,y_1)-{\mathbbm{1}_{\{|z|\leq 1\}}}z\cdot D_x\Gamma(x_1,y_1)\big)}{|z|^{N+2s}}dz\\
			&\qquad+\b C(N,s)\int_{\mathbb{R}^N\setminus B_\delta}\frac{\big(u(x_1+z)-u(x_1)-{\mathbbm{1}_{\{z\leq 1\}}}z\cdot D_x\Gamma(x_1,y_1)\big)}{|z|^{N+2s}}dz\bigg)\\	
			&\qquad+L_1\big(b(x_1)-b(y_1)\big)\cdot \frac{\phi'(x_1-y_1)}{|x_1-y_1|}(x_1-y_1)-b(y_1)\cdot D_y\Gamma(x_1,y_1)+4L_2\|b\|_{\infty,B_1},
		\end{align}
	 This together with \ref{H2} yields using the definition of $F$ that
		\begin{align}
			\frac{{f}(x_1)}{\big|{p}+D_x\Gamma(x_1,y_1)\big|^{\gamma(x_1)}}&\leq\bigg(\a F(Y)+\b C(N,s)\int_{B_\delta}\frac{\big(-\Gamma(x_1,y_1+z)+\Gamma(x_1,y_1)+{\mathbbm{1}_{\{|z|\le 1\}}}z\cdot D_y\Gamma(x_1,y_1)\big)}{|z|^{N+2s}}dz\\
			&\qquad+\b C(N,s)\int_{\mathbb{R}^N\setminus B_\delta}\frac{\big(u(y_1+z)-u(y_1)-{\mathbbm{1}_{\{|z|\leq 1\}}}z\cdot D_y\Gamma(x_1,y_1)\big)}{|z|^{N+2s}}dz\bigg)-b(y_1)\cdot D_y\Gamma(x_1,y_1)\\
			&\qquad+\a  \mathcal{M}_{\lambda,\Lambda}^+(X-Y)-P+4L_2\|b\|_{\infty,B_1}\\
			&\leq\bigg(\a F(Y)+\b C(N,s)\int_{B_\delta}\frac{\big(-\Gamma(x_1,y_1+z)+\Gamma(x_1,y_1)+{\mathbbm{1}_{\{|z|\leq 1\}}}z\cdot D_y\Gamma(x_1,y_1)\big)}{|z|^{N+2s}}dz\\
			&\qquad+\b C(N,s)\int_{\mathbb{R}^N\setminus B_\delta}\frac{\big(u(y_1+z)-u(y_1)-{\mathbbm{1}_{\{|z|\leq 1\}}}z\cdot D_y\Gamma(x_1,y_1)\big)}{|z|^{N+2s}}dz\bigg)-b(y_1)\cdot D_y\Gamma(x_1,y_1)\\
			&\qquad+\a\mathcal{M}_{\lambda,\Lambda}^+(X-Y)-P+4L_2\|b\|_{\infty,B_1}.
		\end{align}
		Further, using the estimate \eqref{in eq 1} in above inequality yields
		\begin{align}
			\frac{{f}(x_1)}{\big|{p}+D_x\Gamma(x_1,y_1)\big|^{\gamma(x_1)}}&\leq \bigg(\a F(Y)+\b C(N,s)\int_{B_\delta}\frac{\big(\Gamma(x_1,y_1+z)-\Gamma(x_1,y_1)+{\mathbbm{1}_{\{|z|\le 1\}}}z\cdot D_y\Gamma(x_1,y_1)\big)}{|z|^{N+2s}}dz\\
			&\qquad+\b C(N,s)\int_{\mathbb{R}^N\setminus B_\delta}\frac{\big(u(y_1+z)-u(y_1)-{\mathbbm{1}_{\{|z|\leq 1\}}}z\cdot D_y\Gamma(x_1,y_1)\big)}{|z|^{N+2s}}dz\bigg)-b(y_1)\cdot D_y\Gamma(x_1,y_1)\\ &\qquad+\alpha\Lambda(N-1)(4L_2+2i)+\alpha\lambda\left(-4(1+\theta)\theta L_12^{\theta-1}\phi_0+4L_2+2i\right)-P+4L_2\|b\|_{\infty,B_1}\\
			&\leq\frac{{f}(y_1)}{\big|{p}-D_y\Gamma(x_1,y_1)\big|^{\gamma(y_1)}}+\alpha\Lambda(N-1)(4L_2+2i)+\alpha\lambda\left(-4(1+\theta)\theta L_12^{\theta-1}\phi_0+4L_2+2i\right)\\
			&\qquad-P+4L_2\|b\|_{\infty,B_1}.
		\end{align}
		Next, for some positive constant $C,$ we have
		\begin{align}
			-Ca^{({\inf_{ B_1}\gamma)}}\varepsilon&\leq -\frac{\varepsilon}{{\big|{p}+D_x\Gamma(x_1,y_1)\big|^{(\inf_{B_1}\gamma)}}}\\
			&\leq-\frac{\varepsilon}{{\big|{p}+D_x\Gamma(x_1,y_1)\big|^{\gamma(x_1)}}}\\
			&\leq -\frac{\|{f}\|_{\infty,B_1}}{{\big|{p}+D_x\Gamma(x_1,y_1)\big|^{\gamma(x_1)}}}\\
			&\leq \frac{{f}(x_1)}{\big|{p}+D_x\Gamma(x_1,y_1)\big|^{\gamma(x_1)}}\\
			&\leq \frac{{f}(y_1)}{\big|{p}-D_y\Gamma(x_1,y_1)\big|^{\gamma(y_1)}}\\
			&\quad+\alpha\Lambda(N-1)(4L_2+2i)+\alpha\lambda\left(-4(1+\theta)\theta L_12^{\theta-1}+4L_2+2i\right)-P+4L_2\|b\|_{\infty,B_1}\\
			&\leq \frac{\varepsilon}{\big|{p}-D_y\Gamma(x_1,y_1)\big|^{\gamma(y_1)}}\\
			&\quad+\alpha\Lambda(N-1)(4L_2+2i)+\alpha\lambda\left(-4(1+\theta)\theta L_12^{\theta-1}\phi_0+4L_2+2i\right)-P+4L_2\|b\|_{\infty,B_1}\\
			&\leq \frac{\varepsilon}{\big|{p}-D_y\Gamma(x_1,y_1)\big|^{(\inf_{B_1}\gamma)}}\\
			&\quad+\alpha\Lambda(N-1)(4L_2+2i)+\alpha\lambda\left(-4(1+\theta)\theta L_12^{\theta-1}\phi_0\phi_0+4L_2+2i\right)-P+4L_2\|b\|_{\infty,B_1}\\
			&\leq Ca^{(\inf_{B_1}\gamma)}\varepsilon+\alpha\Lambda(N-1)(4L_2+2i)+\alpha\lambda\left(-4(1+\theta)\theta L_12^{\theta-1}\phi_0+4L_2+2i\right)-P+4L_2\|b\|_{\infty,B_1},
		\end{align}
	where in the first and last inequalities, we used \eqref{p++} and \eqref{p--}, respectively.
   		This further implies
		\begin{align}
			4\alpha(1+\theta)\theta\lambda L_12^{\theta-1}\phi_0&\leq 2Ca^{(\inf_{B_1}\gamma)}\varepsilon+\alpha\Lambda(N-1)\big(4L_2+2i\big)+\alpha\lambda\big(4L_2+2i\big)-P+4L_2\|b\|_{\infty,B_1}\\
			&=2Ca^{(\inf_{B_1}\gamma)}\varepsilon+\alpha\big(\Lambda(N-1)+\lambda\big)\big(4L_2+2i\big)-P+4L_2\|b\|_{\infty,B_1}\\
			&=2C_1{L_1}^{-(\inf_{B_1}\gamma)}\varepsilon+\alpha\big(\Lambda(N-1)+\lambda\big)\big(4L_2+2i\big)-P+4L_2\|b\|_{\infty,B_1},
		\end{align} 
		for some $C_1>0.$ Now, choosing $L_1$ large enough, depending only on $\lambda,\Lambda,N,$ $L_2$, and using the condition $\inf_{B_1}\gamma\geq 0,$ we get a contradiction. This implies that $u$ is Lipschitz continuous in $B_r,$ for $r\in(0,{1}/{2}).$
		\qed\\
		
		Next, we consider the case when $p\leq\frac{1}{a}.$
		\begin{lem}\label{lemma 4'}
			Let \ref{H3}, \ref{H1} and \ref{H2} hold and $s\in\big({1}/{2},1\big).$ Let $u$ be a viscosity solution of 
			\begin{align}
				|Du+p|^{\gamma(x)}\big(\a F\big(D^2u\big)-\b(-\Delta)^su\big)+b\cdot Du|Du+p|^{\gamma(x)}=f \text{ in }B_1
			\end{align}
			with ${\text{osc}}_{B_1}u\leq 1$ and $\|f\|_{\infty,B_1}\leq\varepsilon$ for some $\varepsilon>0.$ If $|p|\leq 
			{1}/{a},$ 
			then $u$ is Lipschitz continuous in $B_r$ with
			\begin{align}
				[u]_{C^{0,1}(B_{r})}\leq C=C\left(\lambda,\Lambda,N,r,\|f\|_{\infty,B_1},\|u\|_{\infty,B_1}\right).
			\end{align}
		\end{lem} 
		\noindent \textbf{Proof of Lemma \ref{lemma 4'}.}
		The proof has exactly the similar lines of proof as in the previous lemma. In this case, we have $|p|\leq {1}/{a}=5cL_1$ (as fixed in Lemma \ref{lemma 4}). Taking $L_1$ large enough so that 
		\begin{align}
			|p+D_x\Gamma|,\,|p-D_y\Gamma|\geq 1.
		\end{align}
		The rest of the proof follows as before. For the sake of brevity, we omit the details.\qed
		
		\begin{lem}\label{lemma .5}
			For every $M>0,\,\varepsilon,\,\delta\in(0,1)$ and a modulus of continuity $\omega\colon [0,+\infty)\To\mathbb{R}_+,$ there exist $k\in(0,1)$ and $s_0\in({(1+\delta)}/{2},1)$ such that for each $p\in\mathbb{R}^N,$ $s\in (s_0,1)$ and a viscosity solution $u$ to 
			\begin{align}
				-k\leq |p+Du|^{\gamma(x)}\big(\a F\big(D^2u\big)-\b\big(-\Delta\big)^su\big)+b\cdot (Du+p)|p+Du|^{\gamma(x)}\leq k \text{ in }B_1,
			\end{align}
			satisfying\\
			(i) $|u(x)-u(y)|\leq \omega(|x-y|)$ for $x,\,y\in\overline{B}_1,$\\
			(ii) $|u(x)|\leq M(1+|x|^{1+\delta}),\,x\in\mathbb{R}^N,$\\
			and given any $\varepsilon>0$, there exists $\eta(\varepsilon)>0$ such that if
			\begin{align}
				\|k(1+|p|)^{-\gamma(x)}\|_{\infty,B_1}\leq \eta,
			\end{align}
			then
			\begin{align}
				\displaystyle{\sup_{B_\frac{1}{2}}}\,|u-h|\leq \varepsilon,
			\end{align}
			for a viscosity solution $h$ to
			\begin{align}\label{1.10}
				\big(\a F\big(D^2u\big)+\b\Delta u\big)+b(x)\cdot (Du+p)=0 \text{ in }B_1,
			\end{align}
		with $h\in C^{1,\overline{\delta}}(B_{\frac{3}{4}})$ for some $\overline{\delta}\in (0,1)$ 
		\end{lem}
		
		\noindent \textbf{Proof of Lemma \ref{lemma .5}.}
		Let if possible, assume that there exists some $M>0,\varepsilon,\delta\in(0,1)$ a modulus of continuity $\omega,$ and sequence $p_n\To p$ in $\mathbb{R}^N,$ $k_n\To 0,$ $s_n\To 1$ with a sequence of functions $\{u_n\}$ satisfying (i) and (ii) listed above along with
		\begin{align}
			\|k_n(|1+|p_n|)^{-\gamma_n(x)}\|_{\infty,B_1}\leq \frac{1}{n},
		\end{align}
		with $0\leq \gamma_n\leq \sup \gamma.$ Moreover, let for each $n\in \mathbb{N},$ $\{u_n\}$ satisfies the following:
		\begin{align}\label{lim h'}
			-k_n\leq|p_n+Du_n|^{\gamma_n(x)}\big(\a F\big(D^2u_n\big)-\b(-\Delta)^{s_n}u_n\big)+b(x)\cdot (Du_n+p_n)|Du_n+p_n|^{\gamma_n(x)}\leq k_n \text{ in }B_1, 
		\end{align}  
		in the viscosity sense such that
		\begin{align}\label{eq contra}
			\displaystyle{\sup_{B_\frac{1}{2}}}\,|u_n-h|> \varepsilon,
		\end{align}
		for any function $h\in C^{1,\alpha}(B_{\frac{3}{4}}).$

		Note that using (i), we have by Arzel\'a-Ascoli theorem that $u_n$ converges locally uniformly in $B_1$ up to a subsequence to a function, say $u_\infty.$ Also, (ii) assures the uniform integrability condition and we have by Vitali's theorem that one can pass the limit. We need to show that the limiting function $u_\infty$ is a solution of
		\begin{align}\label{eq lim h}
			\big(\a F\big(D^2u\big)+\b\Delta u\big)+b(x)\cdot (Du+p)=0 \text{ in }B_1.
		\end{align}
		
		We closely follow the arguments of Bronzi et al. \cite{Bronzi}, who considered the fully nonlinear local operator with gradient degeneracy, $|\nabla u|^{\gamma(x)}F(D^2u).$ These arguments have been also exploited by Fang et al. \cite{Fang}. In our case, we also have a nonlocal term $(-\Delta)^s$ and a first-order term $b(x)\cdot Du$ present. We re-write \eqref{lim h'} as follows:
		\begin{align}\label{any n}
			-k_n
			\leq &\left(1+|p_n|\right)^{\gamma_n(x)}\left|\frac{p_n}{1+|p_n|}+\frac{1}{1+|p_n|}Du_n\right|^{\gamma_n(x)}\bigg(\a F\big(D^2u_n\big)-\b(-\Delta)^{s_n}u_n+b\cdot \left({Du_n+p_n}\right)\bigg)\leq k_n.
		\end{align}
		It is easy to see that
		\begin{align}
			\xi_n:=\frac{p_n}{1+|p_n|}\To \xi,\,\,\,\,\, \text{ and }\,\,\,\,\, \mu_n:=\frac{1}{1+|p_n|}\To \mu
		\end{align}
		up to a subsequence with $|\xi|\leq 1$ and $\mu\in [0,1].$ Let $A$ be a real symmetric matrix and $v$ be some constant vector. Consider
		\begin{align}
			P(x)\coloneqq\frac{1}{2}\langle A(x-y),(x-y)\rangle+v\cdot(x-y)+u_\infty(y)
		\end{align}
		be a quadratic polynomial which touches $u_\infty$ from below at some point $y$ in $B_{{3}/{4}}.$ Without any loss of generality, let $|y|=u_\infty(y)=0.$ 
We only aim to show that $u_\infty$ is a supersolution to \eqref{eq lim h}. The subsolution part can be proved using analogous arguments. Fix $r\in (0,1)$, we define 
		\begin{align}
			P_j(x)\coloneqq\frac{1}{2}\langle A(x-x_j),(x-x_j)\rangle+v\cdot (x-x_j)+u_j(x_j)
		\end{align}
		with $x_j$ satisfying
		\begin{align}
			P_j(x_j)-u_j(x_j):=\displaystyle{\max_{x\in B_r}}(P_j(x)-u_j(x)).
		\end{align}
		This is possible since $u_n\To u_\infty$ locally uniformly in $B_1.$ Consider
		\begin{align}
			\widetilde{P}_j=\begin{cases}
				P_j &\text{ in }B_{\delta'}(x_j)\\[1mm]
				u_j &\text{ in }\mathbb{R}^N\setminus B_{\delta'}(x_j)		
			\end{cases}
		\end{align}
		We have for the choice $n=j$ in \eqref{any n},
		\begin{align}\label{n=j 12}
			(1+|p_j|)^{\gamma_j(x_j)}\left|\frac{p_j}{1+|p_j|}+\frac{1}{1+|p_j|}v\right|^{\gamma_j(x_j)}\left(\a F\big(A\big)-\b(-\Delta)^{s_j}\widetilde{P}_j(x_j)+b(x_j)\cdot\left({v+p_j}\right)\right)\leq k_j.
        \end{align}
		It gives
		\begin{align}
			\left|\frac{p_j}{1+|p_j|}+\frac{1}{1+|p_j|}v\right|^{\gamma_j(x_j)}\left(\a F\big(A\big)-\b(-\Delta)^{s_j}\widetilde{P}_j(x_j)+b(x_j)\cdot (DP_j(x_j)+p_j)\right)\leq \frac{k_j}{[1+|p_j|]^{\gamma_j(x_j)}}.
		\end{align}
		Let 
		\begin{align}
			\left|\frac{p_j}{1+|p_j|}+\frac{1}{1+|p_j|}v\right|^{\gamma_j(x_j)}\To \widetilde{\beta}\in [0,\infty],
		\end{align}
		up to a subsequence. Consider the first case, i.e., $\widetilde{\beta}>0,$ then taking the limit $j\To\infty$ in \eqref{n=j 12}, we get
		\begin{align}
			\big(\alpha F\big(A\big)+\beta\text{ trace}(A)+b(0)\cdot (DP(0)+p)\big)\leq 0.
		\end{align}
		This further yields
		\begin{align}\label{leq 0}
			\big(\alpha F\big(A\big)+\beta\text{ trace}(A)+b\cdot (Du_\infty+p)\big)\leq 0 \text{ in }B_{\frac{3}{4}}.
		\end{align}
		For the complementary case, $v=0,$ $|p_n|\To 0.$ Consider $$v_n=u_n+\langle p_n,x\rangle.$$ 
		We only need to analyze the case $v=0=p_n.$ Also, we consider only the case when $A$ is not negative semi-definite, since for $A\leq 0,$ we immediately have
		\eqref{leq 0} (since $v=0$). Thus, we only need to consider the case when the invariant space formed by eigenvectors associated with positive eigenvalues, say $E:=\,$span$\,(e_1,e_2,\dots e_n)\in S^{N-1}$ is non-empty. Consider an orthogonal sum, $\mathbb{R}^N=E\oplus G.$ Let $P_Ex$ be denote the orthogonal projection of $x$ on $E$. Define the test function:
		\begin{align}
			\Phi(x)\coloneqq\frac{1}{2}\langle Ax,x\rangle+m\displaystyle{\sup_{e\in S^{N-1}}}\langle P_Ex,e\rangle,
		\end{align}  
		for $m>0$ small enough. Now, since we have that $u_j\To u_\infty$ locally uniformly and $\frac{1}{2}\langle Ax,x\rangle$ touches $u_\infty$ at $0$ from below so we can choose $m$ small enough such that $\Phi$ touches $u_j$ from below at some point $x_j^m\in B_r.$  Moreover, if $x_j^m\in G,$ then
		\begin{align}
			\Phi(x)\coloneqq \frac{1}{2}\langle Ax,x\rangle+m\langle P_Ex,e\rangle
		\end{align}
		touches $u_j$ at $x^m_j$ for any choice of $e.$ We define
		\begin{align}
			\widetilde{\Phi}\coloneqq \begin{cases}
				\Phi &\text{ in }B_{\delta'}(x_j^m)\\[1mm]
				u_j &\text{ in }\mathbb{R}^N\setminus B_{\delta'}(x_j^m).	
			\end{cases}
		\end{align}
		Note that here $\delta'$ is a small number which we can choose accordingly smaller than the previous choice. It is easy to observe that
		\begin{align}
			D\Phi(x_j^m)=Ax_j^m+kP_Ee\,\text{ and }\,D^2\Phi(x_j^m)=A.
		\end{align}
		For the choice $n=j,$ we get
		\begin{align}
			k_j&\geq |D\Phi(x_j^m)|^{\gamma_j(x_j^m)}\big(\a F\big(D^2\Phi(x_j^m)\big)-\b\big(-\Delta\big)^{s_j}\widetilde{\Phi}(x_j^m)\big)+b\cdot D\Phi(x_j^m)|D\Phi(x_j^m)|^{\gamma_j(x_j^m)}\\
			&=|Ax_j^m+me|^{\gamma_j(x_j^m)}\big(\a F\big(D^2\Phi(x_j^m)\big)-\b\big(-\Delta\big)^{s_j}\widetilde{\Phi}(x_j^m)+b\cdot D\Phi(x_j^m)\big).
		\end{align}
		Next, following the arguments of \cite{Bronzi} and taking supremum over all $e\in S^{N-1}$ followed by letting $j\To\infty$ provides
		\begin{align}
			\a F\big(D^2\Phi(x_j^m)\big)+\beta \Delta\Phi(x_j^m) +b\cdot D\Phi(x_j^m)\leq 0.	
		\end{align}  
		Finally, consider the case when $P_E(x_j^m)\neq 0.$ Consider
		\begin{align}
			\rho_j^m=\frac{P_E(x_j^m)}{|P_E(x_j^m)|},
		\end{align} 
		where $$|P_E(x_j^m)|\coloneqq\displaystyle{\sup_{e\in S^{N-1}}\langle P_E(x_j^m),e\rangle}$$ denotes the Euclidean norm of $P_E(x_j^m).$ It is easy to see that $|P_E(x)|$ is a smooth function with
		\begin{align}
			D(|P_E(x)|)_{x={x_j^m}}=\rho_j^m \text{ and } D^2(|P_E(x)|)_{x={x_j^m}}=I-\rho_j^m\otimes\rho_j^m.
		\end{align}
		Hence, we have for the choice $n=j$, following the arguments of \cite{Bronzi} that
		\begin{align}\label{kId}
			\big(\a F(A+m(I-\rho_j^m\otimes \rho_j^m))-\b(-\Delta)^{s_j}\widetilde{\Phi}(x_j^m)\big)+b(x_j^m)\cdot D\Phi(x_j^m)\leq k_j|Ax_j^m+m\rho_j^m|^{-\gamma_j(x_j^m)}.
		\end{align}
		Now, using the fact that $m(I-\rho_j^k\otimes \rho_j^k)$ is positive semi-definite, we get
		\begin{align}
			F(A+m(I-\rho_j^m\otimes \rho_j^m))\geq F(A).
		\end{align}
		Using this in \eqref{kId} yields
		\begin{align}
			\big(\a F(A)-\b(-\Delta)^{s_j}\widetilde{\Phi}(x_j^m)\big)+b(x_j^m)\cdot D\Phi(x_j^m)\leq k_j|Ax_j^m+m\rho_j^m|^{-\gamma_j(x_j^m)}.
		\end{align}
		Finally, letting $j\To\infty,$ concludes
		\begin{align}
			\alpha F(A)+\beta\, \text{trace}(A)+b(0)\cdot D\Phi(0)\leq 0.	
		\end{align} 
		Hence, $u_\infty$ is a viscosity supersolution of
		\begin{align}
			\alpha F\big(D^2u\big)+\beta \Delta u+b\cdot (Du+p)=0.	
		\end{align} 
		Similarly, one may check the subsolution part. Hence, we have that $u_\infty$ is a viscosity solution of \eqref{eq lim h}. By Remark 3.5 \cite{C1beta}, we have that $u_\infty\in C^{1,\overline{\delta}}(B_{{3}/{4}})$ for some $0<\overline{\delta}<1.$ Finally, taking $h=u_\infty,$ we get a contradiction to \eqref{eq contra} for all $n$ large enough. This proves the claim.\qed
		
		
	\begin{lem}\label{Lemma .6}
			Let \ref{H3}, \ref{H1} and \ref{H2} hold. There exist $\tau,\delta,k\in (0,1),$ $s_0\in\big({1}/{2},1\big)$ such that for any bounded viscosity solution $u$ of
			\begin{align}\label{in B1}
				-k\leq|Du|^{\gamma(x)}\left(\a F\big(D^2u\big)-\b(-\Delta)^su\right)+b\cdot Du|Du|^{\gamma(x)}\leq k \text{ in }B_1,
			\end{align} 
			for $s\in (s_0,1)$ with $\|u\|_{\infty,\mathbb{R}^N}\leq 1,$ there exists a sequence of affine functions $t_n=a_n+p_n\cdot x$ with $\{a_n\}\subset\mathbb{R}$ and $\{p_n\}\subset\mathbb{R}^N$ such that
			
			\begin{align}\label{.2}
				\begin{cases}
					\displaystyle{\sup_{B_{\tau^n}}}\,|u-t_n|&\leq \tau^{n(1+\delta)},\\[1mm]
					|a_{n+1}-a_n|&\leq C\tau^{(1+\delta)n},\\[1mm]
					|p_{n+1}-p_n|&\leq C\tau^{\delta n},
				\end{cases}
			\end{align} 
			for some positive constant $C.$ 
		\end{lem}
		\noindent \textbf{Proof of Lemma \ref{Lemma .6}.}
		Let $h\in C(B_1)$ (see Lemma \ref{lemma .5} for such a function) be viscosity solution
		\begin{align}
			\a F\big(D^2u\big)+\beta\Delta u+b\cdot Du=0,	
		\end{align} 
		with $\|h\|_{\infty,B_1}\leq C'$ for some $C'> 0.$ We set
		\begin{align}\label{t_1}
			t_1=h(0)+Dh(0)\cdot x.	
		\end{align} 
		Let us consider universal constants $\overline{A},\,\overline{\delta}\in (0,1)$ such that 
		\begin{align}
			|h(0)|,|Dh(0)|,\|h\|_{C^{1,\overline{\delta}}(B_\tau)}\leq \overline{A}\,\,\,\,\,\, \text{ and }\,\,\,\,\,\, \displaystyle{\sup_{B_\tau}}|h(x)-h(0)-Dh(0)\cdot x|\leq \overline{A}|x|^{1+\overline{\delta}},  
		\end{align}
		for $B_r\Subset B_{\frac{3}{4}}.$ Moreover, since $\|u-h\|_{\infty,B_{{1}/{2}}}\leq \varepsilon$, so we have \begin{align}
			\displaystyle{\sup_{B_r}}\,|u(x)-h(0)-b\cdot Dh(0)|\leq \varepsilon+\overline{A}|x|^{1+\overline{\delta}}.
		\end{align}
		Fix $\delta\in (0,\overline{\delta})\cap\left(0,\frac{1}{1+\|\gamma^+\|_{\infty}}\right)$ such that $1-\delta(1+\sup_{B_1}\gamma)>0,$ and $0<\tau<1$ small enough such that
		\begin{align}\label{A bar}
			\tau^\delta\leq \frac{1}{4},\, \tau^{\overline\delta-{\delta}}(1+\overline{A})\leq \varepsilon_1\,\,\, \text{ and }\,\, (1+2\overline{A})\tau^{1-\delta(1+\|\gamma^+\|_\infty)}\leq 1
		\end{align}
		for a small number $\varepsilon_1>0.$ We take $k$ small enough and $s_0$ close enough to $1$ such that for any solution $u$ of $\eqref{.2}$ satisfying the condition (ii) of Lemma \ref{lemma .5}, we have the existence of a function $h\in C^{1,\overline{\delta}}(B_{3/4})$ that solves \eqref{1.10} in the viscosity sense with $$\|u-h\|_{\infty,B_{{1}/{2}}}\leq\varepsilon.$$ Recall that the condition (i) of Lemma \ref{lemma .5} follows by Theorem \ref{Main 1}.
		Also, set 
		\begin{align}
			\varepsilon\coloneqq\frac{\tau^{1+\delta}}{2}.
		\end{align}
		We prove the result by induction. Consider the auxiliary function:
		\begin{align}\label{w u}
			w_{n+1}(x)\coloneqq\frac{u(\tau^nx)-t_n(\tau^nx)}{\tau^{n(1+\delta)}},
		\end{align}
		with $t_0=0$ and $$t_n(x)=a_n+p_n\cdot x,$$ where $a_n,\,p_n$ to be defined later. It can be re-written as
		\begin{align}\label{u w}
			u(\tau^nx)=\tau^{n(1+\delta)}w_{n+1}(x)+t_n(\tau^nx).
		\end{align}
		We have the following relations in the viscosity sense$\colon$
		\begin{align}
			\tau^nDu(\tau^nx)=\tau^{n(1+\delta)}Dw_{n+1}(x)+\tau^np_n,
		\end{align}
		i.e.,
		\begin{align}\label{grad u w}
			Du(\tau^nx)=\tau^{n\delta}Dw_{n+1}(x)+p_n.
		\end{align}
		Also,
		\begin{align}
			D^2u(\tau^nx)=\tau^{n(\delta-1)}D^2w_{n+1}(x).
		\end{align}
		This grants
		\begin{align}\label{pcci u w}
			F\big(D^2u(\tau^nx)\big)=F\big(\tau^{n(\delta-1)}D^2w_{n+1}(x)\big).
		\end{align}
		On the other hand, since
		\begin{align}
			-(-\Delta)^su(\tau^nx)=\frac{C(N,s)}{2}\int_{\mathbb{R}^N}\frac{u(\tau^nx+y')+u(\tau^nx-y')-2u(\tau^nx)}{|y'|^{N+2s}}dy',
		\end{align}
		so substituting $y'=\tau^ny$ yields
		\begin{align}\label{u taun}
			-(-\Delta)^su(\tau^nx)&=\frac{C(N,s)}{2}\int_{\mathbb{R}^N}\frac{u(\tau^nx+\tau^ny)+u(\tau^nx-\tau^ny)-2u(\tau^nx)}{|\tau^ny|^{N+2s}}\tau^{nN}dy\\
			&=\frac{C(N,s)}{2}\tau^{-2ns}\int_{\mathbb{R}^N}\frac{u(\tau^n(x+y))+u(\tau^n(x-y))-2u(\tau^nx)}{|y|^{N+2s}}dy.
		\end{align}
		Using \eqref{u w} in \eqref{u taun} infers
		\begin{align}\label{3.23}
			-(-\Delta)^su(\tau^nx)&=\frac{C(N,s)}{2}\tau^{n(1+\delta)}\tau^{-2ns}\int_{\mathbb{R}^N}\frac{w_{n+1}(x+y)+w_{n+1}(x-y)-2w_{n+1}(x)}{|y|^{N+2s}}dy\\
			&=-\tau^{n(1+\delta-2s)}(-\Delta)^sw_{n+1}(x).
		\end{align}
		Using \eqref{grad u w}, \eqref{pcci u w} and \eqref{3.23} in \eqref{in B1} yields
		\begin{align}
			-k\leq \big|\tau^{n\delta}Dw_{n+1}(x)+p_n\big|^{\gamma_n(x)}\bigg(&\alpha F\big(\tau^{n(\delta-1)}D^2w_{n+1}(x)\big)-\tau^{n(1+\delta-2s)}\beta(-\Delta)^sw_{n+1}(x)\bigg)\\
			&\qquad+b(\tau^nx)\cdot \big(\tau^{n\delta}Dw_{n+1}(x)+p_n\big)\big|\tau^{n\delta}Dw_{n+1}(x)+p_n\big|^{\gamma_n(x)}\leq k \text{ in }B_{\tau^{-n}},
		\end{align}
		where $$\gamma_n(x)=\gamma{(\tau^n(x))}.$$ In other words,
		{\begin{align}\label{w solves}
				\quad-k\leq
				\tau^{n\delta\gamma_n(x)}\big|Dw_{n+1}(x)+\tau^{-n\delta}p_n\big|^{\gamma_n(x)}\bigg(&\alpha F\big(\tau^{n(\delta-1)}D^2w_{n+1}(x)\big)-\tau^{n(1+\delta-2s)}\beta(-\Delta)^sw_{n+1}(x)\\
				&\qquad\qquad\qquad\quad+b(\tau^nx)\cdot \big(\tau^{n\delta}Dw_{n+1}(x)+p_n\big)\bigg)\leq k \text{ in }B_{\tau^{-n}}.
			\end{align}
			It deduces that $w_{n+1}$ given by \eqref{w u} solves \eqref{w solves} in $B_{\tau^{-n}}.$} Now, multiplying $\tau^{n(1-\delta-\delta\gamma_n(x))}$ in both sides of \eqref{w solves}, we get 
		\begin{align}
			-\tau^{n(1-\delta-\delta\gamma_n(x))}k&\leq \big|Dw_{n+1}(x)+\tau^{-n\delta}p_n\big|^{\gamma_n(x)}\bigg(\tau^{n(1-\delta)}\alpha F\big(\tau^{n(\delta-1)}D^2w_{n+1}(x)\big)-\tau^{n(1-\delta)}\tau^{n(1+\delta-2s)}\beta(-\Delta)^sw_{n+1}(x)\\
			&\qquad\quad\qquad\qquad\quad\qquad\qquad\quad\qquad\qquad\qquad\quad\qquad\qquad\quad\qquad+\tau^{n(1-\delta)}b(\tau^nx)\cdot \big(\tau^{n\delta}Dw_{n+1}(x)+p_n\big)\bigg)\\
			&=\big|Dw_{n+1}(x)+\tau^{-n\delta}p_n\big|^{\gamma_n(x)}\bigg(\tau^{n(1-\delta)}\tau^{n(\delta-1)}\alpha F\big(D^2w_{n+1}(x)\big)-\tau^{n(2-2s)}\beta(-\Delta)^sw_{n+1}(x)\\
			& \qquad\quad\qquad\quad\qquad\qquad\qquad\qquad\quad\qquad\qquad\quad\qquad\qquad\qquad\quad+\tau^{n(1-\delta)}b(\tau^nx)\cdot \big(\tau^{n\delta}Dw_{n+1}(x)+p_n\big)\bigg)\\
			&=\big|Dw_{n+1}(x)+\tau^{-n\delta}p_n\big|^{\gamma_n(x)}\bigg(\alpha F\big(D^2w_{n+1}(x)\big)-\tau^{n(2-2s)}\beta(-\Delta)^sw_{n+1}(x)\\
			&\qquad\quad\qquad\quad\qquad\qquad\quad\qquad\qquad\qquad\qquad\quad\qquad\qquad\qquad\quad+\tau^{n(1-\delta)}b(\tau^nx)\cdot \big(\tau^{n\delta}Dw_{n+1}(x)+p_n\big)\bigg)\\
			&\leq \tau^{n(1-\delta-\delta\gamma_n(x))}k, \text{ for } x\in B_{\tau^{-n}}.&
		\end{align}
		This further accords
		\begin{align}\label{More w}
			\\\quad-k	&\leq -\tau^{(1-\delta(1+\sup_{ B_{1}}\gamma(x)))}k\\
			&=-\tau^{(1-\delta(1+\sup_{ B_{\tau^{-n}}}\gamma_n(x)))}k\\
			&\leq -\tau^{(1-\delta(1+\gamma_n(x)))}k\\&\leq \big|Dw_{n+1}(x)+\tau^{-n\delta}p_n\big|^{\gamma_n(x)}\bigg(\alpha F\big(D^2w_{n+1}(x)\big)-\tau^{n(2-2s)}\beta(-\Delta)^sw_{n+1}(x)\\
			&\qquad\qquad\qquad\qquad\qquad\qquad\qquad\qquad\qquad+\tau^{n(1-\delta)}b(\tau^nx)\cdot \big(\tau^{n\delta}Dw_{n+1}(x)+p_n\big)\bigg)\\
			&\leq \tau^{(1-\delta(1+\gamma_n(x)))}k\\
			&\leq \tau^{(1-\delta(1+\sup_{ B_{\tau^{-n}}}\gamma_n(x)))}k\\
			&=\tau^{(1-\delta(1+\sup_{ B_1}\gamma(x)))}k\\
			&\leq k\,\,\,\,\,\, \text{ for } x\in B_{\tau^{-n}}.
		\end{align}
		It implies that $w_{n+1}$ satisfies the inequality \eqref{More w} in the viscosity sense. Next, our aim is to define $t_n$ in such a way that $w_{n+1}$ satisfies the conditions of Lemma \ref{lemma .5}. This would give us the existence of a viscosity solution of \eqref{1.10}, say $\widetilde{h}\in C^{1,\overline{\delta}}(B_{{3}/{4}})$ such that $$\sup_{B_{{1}/{2}}}|w_{n+1}-\widetilde{h}|\leq {\tau^{1+\delta}}/{2}.$$ We follow the arguments of \cite{Topp}. Similar arguments have also been used, for instance, in \cite{Bronzi, Filippis, Huaroto, Jesus}. We define 
		\begin{align}
			\widetilde{t}=\widetilde{h}(0)+D\widetilde{h}(0)\cdot x
		\end{align} 
		and 
		\begin{align}\label{eq t tau}
			{t}_{n+1}(x)=t_n(x)+\tau^{n(1+\delta)}\widetilde{t}(\tau^{-n}x).
		\end{align} 
		In other words, \begin{align}
			{a}_{n+1}=a_n+\tau^{n(1+\delta)}\widetilde{h}(0),
		\end{align}
	and
		 \begin{align}
		 	{p}_{n+1}=p_n+\tau^{n\delta}D\widetilde{h}(0)\cdot x.
		 \end{align} 
		Taking $w_0=u$ and $t_0=0,$ we have the sequences $\{w_n\}$ (defined by \eqref{w u}) and $\{t_n\}$(defined by \eqref{eq t tau}). Let the conditions of Lemma \ref{lemma .5}, in particular, (ii) holds, i.e.,
		\begin{align}\label{grth}
			|u(x)|\leq (1+|x|^{1+\overline{\delta}}),\,x\in\mathbb{R}^N,
		\end{align}  
		 for $\{w_i\}$ for $i=1,2,\dots n.$ We prove the same for $i=n+1.$ We first re-write $w_{n+1}$ from \eqref{w u} in terms of $\widetilde{t}$ and $w_n.$
		\begin{align}\label{w ow}
			w_{n+1}(x)&=\frac{u(\tau^nx)-t_n(\tau^nx)}{\tau^{n(1+\delta)}}\\
			&=\frac{u(\tau^nx)-\big(t_{n-1}(\tau^nx)+{\tau^{(n-1)(1+\delta)}}\widetilde{t}{(\tau^{-(n-1)}\tau^nx)}\big)}{\tau^{n(1+\delta)}}\\
			&=\frac{u(\tau^nx)-\big(t_{n-1}(\tau^nx)+{\tau^{{(n-1)}(1+\delta)}}\widetilde{t}{(\tau x)}\big)}{\tau^{n(1+\delta)}}\\
			&=\frac{\frac{u(\tau^{n-1}\tau x)-t_{n-1}(\tau^{n-1}\tau x)}{\tau^{(n-1)(1+\delta)}}-\widetilde{t}\big(\tau x\big)}{\tau^{(1+\delta)}}\\
			&=\frac{w_n(\tau x)-\widetilde{t}{(\tau x)}}{\tau^{1+\delta}}.
		\end{align}
		Now, if $\tau x\in \overline{B^c}_{{1}/{2}},$ then by \eqref{w ow}, we have
		\begin{align}
			|w_{n+1}(x)|&\leq \frac{\big|w_n(\tau x)\big|+\big|\widetilde{t}{(\tau x)}\big|}{\tau^{1+\delta}}\\
			&=\frac{\big(1+|\tau x|^{1+\overline{\delta}}\big)+\overline{A}(1+\tau|x|)}{\tau^{1+\delta}}\\
			&\leq\frac{(2|x|\tau+|\tau x|^{1+\overline{\delta}})+\overline{A}(3\tau|x|)}{\tau^{1+\delta}}\\
			&\leq \frac{(2^{1+\overline{\delta}}|x|^{1+\overline{\delta}}\tau^{1+\overline{\delta}}+\tau^{1+\overline{\delta}}| x|^{1+\overline{\delta}})+\overline{A}(3^{1+\overline{\delta}}\tau^{1+\overline{\delta}}|x|^{1+\overline{\delta}})}{\tau^{1+\delta}}\\
			&\leq \frac{5|x|^{1+\overline{\delta}}\tau^{1+\overline{\delta}}+9\overline{A}|x|^{1+\overline{\delta}}\tau^{1+\overline{\delta}}}{\tau^{1+\delta}}\\
			&\leq \big(5+9\overline{A}\big)\tau^{\overline{\delta}-\delta}|x|^{1+\overline{\delta}}.
		\end{align}
		Further, using \eqref{A bar}, we get
		$|w_{k+1}|\leq |x|^{1+\overline{\delta}}.$
		In the case when $\tau x\in B_{\frac{1}{2}},$ we have
		\begin{align}
			|w_{n+1}(x)|&\leq \frac{\big|w_n(\tau x)\big|+\big|\widetilde{t}{(\tau x)}\big|}{\tau^{1+\delta}}\\
			&\leq \frac{\big|w_n(\tau x)-\widetilde{h}(\tau x)\big|+\big|\widetilde{h}(\tau x)-\widetilde{t}{(\tau x)}\big|}{\tau^{1+\delta}}\\
			&\leq \frac{\frac{\tau^{1+\delta}}{2}+\overline{A}|\tau x|^{1+\overline{\delta}}}{\tau^{1+\delta}}\text{ (using Lemma \ref{lemma .5})}\\
			&\leq \frac{1}{2}+\overline{A}|\tau |^{\overline{\delta}-\delta}|x|^{1+\overline{\delta}}\\
			&\leq 1+|x|^{1+\overline{\delta}} \text{ (using \eqref{A bar})}.
		\end{align}
		Thus, we have that $w_{n+1}$ satisfies \eqref{grth} for each $n\in\mathbb{N}.$ Moreover, notice that
		\begin{align}
			|p_n|&\leq|p_1|+\sum_{i=2}^n|p_i-p_{i-1}|\\
			&\leq |p_1|+\sum_{i=2}^n Cr^{(i-1)\delta}\\
			&\leq C+\frac{C\tau^\gamma}{1-\tau^\gamma} \text{ (using \ref{t_1})}\\
			&\leq \frac{4}{3}C.
		\end{align}
		We choose $C=\overline{A}.$ Thanks to \eqref{A bar} and (pp. 21 \cite{Bronzi}), we have
		\begin{align}
			\big\|k(1+|r^{-n\delta}p_n|)^{-\gamma_n(x)}\big\|_{\infty}\leq k.
		\end{align}
		Taking $k$ small enough, we have by Lemma \ref{lemma .5} that there exists a function $\widetilde{h}$ for each $w_{n+1}$ such that $\widetilde{h}$ solves \eqref{1.10} and $$\sup_{B_{\tau}}|w_{n+1}-\widetilde{h}|\leq \frac{\tau^{(1+\delta)}}{2}.$$ Also, we see that
		\begin{align}\label{w tilde t}
			\sup_{x\in B_\tau}|w_{n+1}-\overline{t}|&\leq \sup_{x\in B_\tau}|w_{n+1}-\widetilde{h}|+\sup_{x\in B_\tau}|\widetilde{t}-\widetilde{h}|\\
			&\leq \frac{\tau^{(1+\delta)}}{2}+Cr^{1+\overline{\delta}}\\
			&\leq \tau^{(1+\delta)},
		\end{align}
		using $\eqref{A bar}$ for $C=\overline{A}.$
		Further, we notice that 
		\begin{align}
			\sup_{B_{\tau^{n+1}}}|u(x)-t_{n+1}|&=\sup_{B_{\tau^{n+1}}}|\tau^{n(1+\delta)}w_{n+1}(\tau^{-n}x)+t_n(x)-t_{n+1}(x)|\text{ (using \eqref{w u})}\\
			&=\sup_{B_{\tau}}|\tau^{n(1+\delta)}w_{n+1}(x)+t_n(\tau^n x)-t^{n+1}(\tau^nx)|\\
			&=\sup_{B_{\tau}}|\tau^{n(1+\delta)}w_{n+1}(x)-\tau^{n(1+\delta)}\widetilde{t}(x)|\text{ (using \eqref{eq t tau})}\\
			&=\tau^{n(1+\delta)}\sup_{B_{\tau}}|w_{n+1}(x)-\widetilde{t}(x)|\\
			&\leq \tau^{(n+1)(1+\delta)},
		\end{align}
		where in the last line, we used \eqref{w tilde t}. Also, we see that 
		\begin{align}
			|a_{n+1}-a_n|\leq \tau^{n(1+\delta)}|\widetilde{h}(0)|\leq C\tau^{n(1+\delta)} \text{ and }|b_{n+1}-b_n|\leq \tau^{n\delta}|D\widetilde{h}(0)|\leq C\tau^{n\delta}.
		\end{align}
		Hence the claim follows.\qed\\
		
		\noindent \textbf{Proof of Theorem \ref{main 2}.} 
		Consider a function
		\begin{align}\label{u tilde u}
			\widetilde{u}:=Ku=\frac{u}{\|u\|_{\infty,B_1}+\left(\nu^{-1}\|f\|_{\infty,B_1}^{\frac{2s-1}{1+\inf_{B_1}\gamma}}\right)},
		\end{align}
		where $K$ is a positive constant and $\nu\in(0,1)$ to be chosen later. Now, if $K\geq 1,$ we immediately have
		\begin{align}
			\|u\|_{\infty,B_1}\leq 1.
		\end{align}
		and
		\begin{align}
			\nu^{-1}\|f\|^{\frac{2s-1}{1+\inf_{B_1}\gamma}}_{\infty,B_1}\leq 1.
		\end{align}
		This further infers
		\begin{align}
			\|f\|_{\infty,B_1}\leq \nu<1.
		\end{align}
		It implies that one can use Lemma \ref{Lemma .6} for $u.$ Hence, we may assume that $K\leq 1.$ By \eqref{eq 0.1}, we have
		\begin{align}
			\|f\|_{\infty,B_1}\geq \frac{|D\widetilde{u}|^{\gamma(x)}}{\left(\|u\|_{\infty,B_1}+\left(\nu^{-1}\|f\|^{\frac{2s-1}{1+\inf_{B_1}\gamma}}_{\infty,B_1}\right)\right)^{-(1+\gamma(x))}}\left(\a F\big(D^2\widetilde{u}\big)-\b(-\Delta)^s\widetilde{u}+b(x)\cdot D\widetilde{u}\right)\geq -\|f\|_{\infty,B_1}.
		\end{align}
		Note that for $\nu\in(0,1)$ small enough, we get $K\ll1.$
		It further implies
		\begin{align}
			1\gg K^{1+\inf_{B_1}\gamma(x)}\|f\|_{\infty,B_1} 
			&\geq K^{1+\gamma(x)}\|f\|_{\infty,B_1}\\
			&\geq {|D\widetilde{u}|^{\gamma(x)}}\left(\a F\big(D^2\widetilde{u}\big)-\b(-\Delta)^s\widetilde{u}\right)+b(x)\cdot D\widetilde{u}|D\widetilde{u}|^{\gamma(x)}\\
			&\geq -K^{1+\gamma(x)}\|f\|_{\infty,B_1}\\
			&\geq -K^{1+\inf_{B_1}\gamma(x)}\|f\|_{\infty,B_1}\\
			&\gg -1.	
		\end{align}
		Hence, we have that $\widetilde{u}$ satisfies the conditions of Lemma \ref{Lemma .6}. It is immediate to see from \eqref{.2} that $\{a_n\}$ and $\{p_n\}$ are Cauchy sequences and hence convergent. Let $a^*$ and $p^*$ be the limits of $\{a_n\}$ and $\{p_n\},$ respectively. Also, consider $$t^*=a^*+b^*\cdot x.$$ It is immediate that
		\begin{align}
			|t^*-t_n|&\leq \sum_{i=n}^\infty|t_{i+1}-t_i|\\
			&\leq \sum_{i=n}^\infty|a_{i+1}-a_i|+\sum_{i=n}^\infty|p_{i+1}-p_i||x|\\
			&\leq \sum_{i=n}^\infty C\tau^{i(1+\delta)}+\sum_{i=n}^\infty C\tau^{i\delta}|x|.
		\end{align}
		Now, letting $$\tau^{n+1}\leq |x|<\tau^{n}, \text{ for some } n\in\mathbb{N},$$ infers
		\begin{align}
			|\widetilde{u}-t^*|&\leq |\widetilde{u}-t_n|+|t^*-t_n|\\
			&\leq C_0\tau^{n(1+\delta)}+\sum_{i=n}^\infty C\tau^{i(1+\delta)}+\sum_{i=n}^\infty C\tau^{i\delta}|x|\text{ (using \eqref{.2})}\\
			&\leq C_0\tau^{n(1+\delta)}+\sum_{i=n}^\infty C\tau^{i(1+\delta)}+\sum_{i=n}^\infty C\tau^{i(1+\delta)}\\
			&\leq C_0\tau^{n(1+\delta)}+ C_1\tau^{n(1+\delta)}+ C_1\tau^{n(1+\delta)}\\
			&= C_2\tau^{n(1+\delta)}\\
			&=\frac{C_2}{\tau^{1+\delta}}\tau^{(n+1)(1+\delta)}\\
			&\leq {C^*}|x|^{1+\delta},
		\end{align}
		where $C_1,$ $C_2\,(=2C_1+C_0)$ and $C^*=\frac{C_2}{\tau^{1+\delta}}$ are positive constants. This implies
		\begin{align}
			[\widetilde{u}]_{C^{1,\delta}(B_{{1}/{2}})}\leq C^*.
		\end{align}
		This together with \eqref{u tilde u} concludes
		\begin{align}
			[u]_{C^{1,\delta}(B_{{1}/{2}})}&\leq C^*\left(\|u\|_{\infty,B_1}+\left(\nu^{-1}\|f\|_\infty^{\frac{2s-1}{1+\inf_{B_1}\gamma}}\right)\right)\\
			&\leq \widetilde{C}\left(\|u\|_{\infty,B_1}+\left(\|f\|_\infty^{\frac{2s-1}{1+\inf_{B_1}\gamma}}\right)\right),
		\end{align} 
		where $\widetilde{C}$ is a positive constant depending on $\nu$ appearing in \eqref{u tilde u}.\qed

		\section{Funding and/or Conflicts of interests/Competing interests}
		The research of Priyank Oza was financially supported by Council of Scientific $\&$ Industrial Research (CSIR) under the grant no. 09/1031(0005)/2019--EMR--I. The second author thanks DST/SERB for the financial support under the grant CRG/2020/000041.\\
		
		There are no conflict of interests of any type. This manuscript does not use any kind of data.

	\end{document}